\documentclass[11pt]{article}
\usepackage{fullpage}

\usepackage[utf8]{inputenc} % allow utf-8 input
\usepackage[T1]{fontenc}    % use 8-bit T1 fonts
\usepackage{hyperref}       % hyperlinks
\usepackage{url}            % simple URL typesetting
\usepackage{booktabs}       % professional-quality tables
\usepackage{amsfonts}       % blackboard math symbols
\usepackage{nicefrac}       % compact symbols for 1/2, etc.
\usepackage{microtype}      % microtypography
\usepackage{tcolorbox}
\usepackage{diagbox}

\usepackage{enumerate}
\usepackage[shortlabels]{enumitem}
\usepackage{algorithmic}
\usepackage{algorithm}
\usepackage{subfigure}
\usepackage{graphicx} % more modern
\usepackage{caption}
\usepackage{amsmath}
\usepackage{amsthm}
\usepackage{amssymb}
\usepackage{tikz}
\usepackage{tablefootnote}
\usepackage{multirow}
\usepackage{enumerate}
\usepackage{color}
\usepackage{xcolor}
\usepackage{natbib}

\usetikzlibrary{arrows}

\allowdisplaybreaks[4]

\usepackage{mathrsfs}

\usepackage{hyperref}
\usepackage{bm,todonotes}

\allowdisplaybreaks

\hypersetup{
	colorlinks=true,
	filecolor=blue,
	citecolor = blue,
	urlcolor=cyan,
}

\newcommand{\bsig}{\bar{\sigma}}
\newcommand{\cD}{\mathcal{D}}
\newcommand{\cO}{\mathcal{O}}

\newcommand{\EE}{\mathbb{E}}
\newcommand{\RR}{\mathbb{R}}

\newcommand{\xb}{\mathbf{x}}

\newcommand{\bs}{\mathbf{s}}

\newcommand{\qb}{\mathbf{q}}

\newcommand{\cF}{{\mathcal{F}}}

\newcommand{\cC}{{\mathcal{E}}}

\newcommand{\bPi}{\mathbf{\Pi}}
\newcommand{\bbx}{\overline{x}}
\newcommand{\bbs}{\overline{s}}
\newcommand{\bbg}{\overline{g}}

\newcommand{\xbti}{\xb_t^{(i)}}

\newcommand{\qbti}{\qb_t^{(i)}}

\newcommand{\ssgt}{\texttt{SS\_GT}}
\newcommand{\tbeta}{\tilde{\beta}}

\newcommand{\txb}{\widetilde{\xb}}

\newcommand{\tsb}{\widetilde{\bs}}
\newcommand{\tW}{\widetilde{W}}
\newcommand{\tPi}{\widetilde{\bPi}}

\newcommand{\algname}{\texttt{SS\_DSGT}}
\newcommand{\algnamea}{\texttt{ASS\_DSGT}}

\newcommand{\dsgt}{\texttt{DSGT}}
\newcommand{\sgd}{\texttt{SGD}}

\newcommand{\norm}[1]{\left\|#1\right\|}
\newcommand{\dotprod}[1]{\left\langle #1\right\rangle}

\newtheorem{lemma}{Lemma}
\newtheorem{theorem}{Theorem}

\newtheorem{definition}{Definition}
\newtheorem{corollary}[theorem]{Corollary}
\newtheorem{remark}{Remark}
\ifx\assumption\undefined
\newtheorem{assumption}{Assumption}

\title{
Snap-Shot Decentralized Stochastic Gradient Tracking Methods
}

\author{
Haishan Ye\thanks{Center for Intelligent Decision-Making and Machine Learning, School of Management, Xi’an Jiaotong University; email: \texttt{yehaishan@xjtu.edu.cn}. }  \\
	\and
	Xiangyu Chang \thanks{Center for Intelligent Decision-Making and Machine Learning, School of Management, Xi’an Jiaotong University; email: \texttt{xiangyuchang@xjtu.edu.cn}. } 
}
\begin{document}

\maketitle

\begin{abstract}%
In decentralized optimization, $m$ agents form a network and only communicate with their neighbors, which gives advantages in data ownership, privacy, and scalability.
At the same time, decentralized stochastic gradient descent (\texttt{SGD}) methods, as popular decentralized algorithms for training large-scale machine learning models, have shown their superiority over centralized counterparts.
Distributed stochastic gradient tracking~(\dsgt)~\citep{pu2021distributed} has been recognized as the popular and state-of-the-art decentralized \texttt{SGD} method due to its proper theoretical guarantees.
However, the theoretical analysis of \dsgt~\citep{koloskova2021improved} shows that
its iteration complexity is $\tilde{\cO} \left(\frac{\bsig^2}{m\mu \varepsilon} + \frac{\sqrt{L}\bsig}{\mu(1 - \lambda_2(W))^{1/2} C_W \sqrt{\varepsilon} }\right)$, where $W$ is a double stochastic mixing matrix that presents the network topology and $ C_W $ is a parameter that depends on $W$.
Thus, it indicates that the convergence property of \dsgt~is heavily affected by the topology of the communication network.
To overcome the weakness of \dsgt, we resort to the snap-shot gradient tracking skill and propose two novel algorithms.
We further justify that the proposed two algorithms are more robust to the topology of communication networks under similar algorithmic structures and the same communication strategy to \dsgt~.
Compared with \dsgt, their iteration complexity are $\cO\left( \frac{\bsig^2}{m\mu\varepsilon} + \frac{\sqrt{L}\bsig}{\mu (1 - \lambda_2(W))\sqrt{\varepsilon}} \right)$ and $\cO\left( \frac{\bsig^2}{m\mu \varepsilon} + \frac{\sqrt{L}\bsig}{\mu (1 - \lambda_2(W))^{1/2}\sqrt{\varepsilon}} \right)$ which reduce the impact on network topology (no $C_W$). 
\end{abstract}

\section{Introduction}\label{sec:intro}

In this paper, we consider the decentralized optimization problem, where there are $ m $ agents to cooperatively minimize a common objective function $ f(x) $ with the following formulation:
\begin{equation}\label{eq:prob}
	f(x) := \frac{1}{m}\sum_{i=1}^m f_i(x), \quad\mbox{with}\quad f_i(x) := \EE_{\xi^{(i)} \sim \cD^{(i)} } f_i(x, \xi^{(i)}).
\end{equation}
The formulation assumes that the objective function  $ f(x) $ is composed of $ m $-local functions $ f_i(x), i\in[m]=\{1,\dots,m\}$.
The $i$-th agent maintains the private data set $\cD^{(i)}$ and its objective function $f_i(x)$.
The $m$ agents form a connected network and can only communicate with their neighbors.
\citet{shi2015extra,Ye2020,scaman2019optimal} indicate that the decentralized optimization has advantages over traditional centralized optimization in data ownership, privacy, and scalability \citep{nedic2020distributed,kairouz2021advances}.
Especially, \citet{lian2017can} provides the first theoretical analysis that indicates
 decentralized algorithms might outperform centralized algorithms of distributed stochastic gradient descent (\sgd).

%At the same time, the stochastic gradient descent (\sgd) has been the most popular work-horse of training machine learning models~\citep{ruder2016overview}. 
Due to the imminent need to train large-scale machine models, decentralized \sgd~methods are attracting significant attention recently because
they are easy to implement, and the computation cost of each iteration is cheap~\citep{Ye2020,shi2015extra,qu2017harnessing,alghunaim2020decentralized}. 
However, the performance of decentralized \sgd~suffers from the data heterogeneity\citep{lian2017can,koloskova2020unified}, that is, training data is in a non-IID fashion distributed over agents. 
%For example, \citet{lian2017can} provides the first theoretical analysis that indicates\xiangyu{what do you mean a regime???}
 %decentralized algorithms might outperform centralized algorithms for distributed \sgd.
{

%It has been shown that data heterogeneity affects the convergence rate of \texttt{DSGD} \citep{lian2017can,koloskova2020unified}.
}

Recently, the gradient tracking method developed by \citet{di2016next} and \citet{nedic2017achieving} has been widely used to overcome the data heterogeneity challenge.
Many decentralized algorithms based on the gradient tracking method have been proposed \citep{qu2017harnessing,song2021optimal,Ye2020}.
For instance, \citet{pu2021distributed} applied the gradient tracking to the decentralized \sgd~ and proposed the distributed stochastic gradient tracking method (\dsgt). 
\dsgt~effectively conquers the dilemma of data heterogeneity, and its dominant computation complexity is the same as its centralized counterpart. 
Furthermore, \dsgt~also has an advantage over centralized \sgd~in communication complexity.

However, the performance of \dsgt~is heavily affected by the topology of the communication network through which  the agents exchange information.
For the $ L $-smooth and $ \mu $-strongly convex functions, 
\dsgt~has the following iteration complexity~\citep{pu2021distributed} to achieve $ \varepsilon $-suboptimality
\begin{equation}\label{eq:com1}
	\tilde{\cO} \left(\frac{\bsig^2}{m \mu\varepsilon} + \frac{\sqrt{L}\bsig}{\mu(1 - \lambda_2(W))^{3/2} \sqrt{\varepsilon} }\right),
\end{equation}
where $ \bsig^2 $ is the upper bound on the variance of the stochastic  noise (see Assumption~\ref{ass:noise}) and  $ \lambda_2(W) $ is the second largest eigenvalue of the double stochastic mixing matrix $ W $.
The above equation shows that when $ \lambda_2(W) $ is close to one, \dsgt~still suffers from poor performance.
Recently, \citet{koloskova2021improved} improved the convergence analysis of \dsgt, and obtain the following complexity
\begin{equation}\label{eq:com2}
		\tilde{\cO} \left(\frac{\bsig^2}{m\mu \varepsilon} + \frac{\sqrt{L}\bsig}{\mu(1 - \lambda_2(W))^{1/2} C_W \sqrt{\varepsilon} }\right),
\end{equation} 
where $ C_W $ is a parameter no smaller than $ 1 - \lambda_2(W) $.
\citet{koloskova2021improved} showed that for a large number of communication networks, $ C_W $ is a constant independent of $ \lambda_2(W) $.
In these cases, Eq.~\eqref{eq:com2} provides a better complexity than Eq.~\eqref{eq:com1}.
Unfortunately, in the general case, $ C_W $ is no longer a constant. 
Eq.~\eqref{eq:com2} may even reduce to Eq.~\eqref{eq:com1} in the worst cases.
Thus, the result in Eq.~\eqref{eq:com1} is the best iteration complexity of \dsgt~for general cases.
It is still an open question:
\textit{can \dsgt~achieve lower communication and computation complexities than Eq.~\eqref{eq:com1} for all communication networks? }
\citet{koloskova2021improved} also proposed an open problem: \textit{is the parameter $ C_W $ in Eq.~\eqref{eq:com2} tight in general for \dsgt?}

Instead of answering the above open questions, we design two novel decentralized stochastic gradient descent tracking algorithms in this paper.
We will justify that the proposed algorithm without extra inner communication loops can achieve lower complexities than Eq.~\eqref{eq:com1}, which take the same communication strategy as \dsgt.
We first extend the ``snap-shot'' gradient tracking method proposed by \citet{song2021optimal} to the stochastic 
 gradient descent. 
 Then we propose a snap-shot decentralized stochastic gradient tracking~(\algname) algorithm accordingly.
\algname~is shown that has the following iteration complexity
\begin{equation*}
	\cO\left( \frac{\bsig^2}{m\mu\varepsilon} + \frac{\sqrt{L}\bsig}{\mu (1 - \lambda_2(W))\sqrt{\varepsilon}} \right),
\end{equation*}
which is better than the one shown in Eq.~\eqref{eq:com1}.
In addition, we use the loopless Chebyshev acceleration \citep{song2021optimal} to improve the performance of \algname~and further we propose \algnamea~which has a better complexity 
\begin{equation} \label{eq:com3}
	\cO\left( \frac{\bsig^2}{m\mu \varepsilon} + \frac{\sqrt{L}\bsig}{\mu (1 - \lambda_2(W))^{1/2}\sqrt{\varepsilon}} \right).
\end{equation}
To the best of our knowledge, our algorithms achieve the best iteration complexities for the decentralized \sgd~without inner communication loops. 
Our result in Eq.~\eqref{eq:com3} shows that one can design a decentralized \sgd~method similar to \dsgt~but without dependency on the parameter $ C_W $.

\section{Notation and Assumptions}

Let $\xb$ and $\bs$ be two $m\times d$ matrices whose $i$-th rows $\xb^{(i)}$ and $\bs^{(i)}$ are the local copy of the decision and gradient-tracking variables  for the $i$-th agent, respectively. 
Accordingly, we define the averaging variables 
\begin{equation}
	\bbx := \frac{1}{m}\sum_{i=1}^m \xb^{(i)} = \frac{1}{m} \mathbf{1}^\top \xb \in\RR^{1\times d}, \quad \bbs := \frac{1}{m}\mathbf{1}^\top \bs \in\RR^{1\times d},
\end{equation}
where $\mathbf{1}$ denotes the vector with all entries equal to $1$. 
Now we introduce the projection matrix
\begin{equation}\label{eq:Pi_def}
	\bPi = \mathbf{I}_m - \frac{\mathbf{1}\mathbf{1}^\top}{m}.
\end{equation}
Using the projection matrix $ \bPi $, we can represent that
\begin{equation*}
	\norm{\xb - \mathbf{1}\bbx} = \norm{ \xb - \frac{\mathbf{1}\mathbf{1}^\top}{m} \xb  } = \norm{\bPi \xb}, \quad\mbox{and} \quad \norm{\bs - \mathbf{1}\bbs} = \norm{\bPi\bs}.
\end{equation*}
We denote an aggregate objective function:
\begin{equation}
	F(\xb) := \sum_{i=1}^m f_i(\xb^{(i)})
\end{equation}
and its aggregate gradient
\begin{equation*}
	\nabla F(\xb) := [\nabla f_1(\xb^{(1)}), \nabla f_2(\xb^{(2)}),\dots, \nabla f_m(\xb^{(m)})]^\top \in\RR^{m\times d}.
\end{equation*}
In addition, let $\xi := [\xi^{(1)}, \xi^{(2)}, \dots, \xi^{(m)}] \in \RR^m$ and 
\begin{equation*}
    \nabla F(\xb, \xi) := \left[\nabla f_1(\xb^{(1)}, \xi^{(1)}), \nabla f_2(\xb^{(2)}, \xi^{(2)}), \dots, \nabla f_m(\xb^{(m)}, \xi^{(m)}) \right].
\end{equation*}

Throughout this paper, we use $ \norm{\cdot} $ to denote the ``Frobenius'' norm. 
That is, for a matrix $ \xb \in \RR^{m\times d} $, it  holds that
\begin{align*}
	\norm{\xb}^2 = \sum_{i=1, j=1}^{m, d} \left(\xb^{(i, j)}\right)^2.
\end{align*}
Furthermore, we use $ \norm{\xb}_2 $ to denote the spectral norm which is the largest singular value of $ \xb $.
For vectors $ x, y \in \RR^d $, we use $ \dotprod{x, y} $ to denote the standard inner product of $ x $ and $ y $.

Now we introduce several assumptions that will be used throughout this paper.
First, we state an assumption that the stochastic gradients have bounded noise.
\begin{assumption}[Bounded Noise]\label{ass:noise}
	We assume that there exists constant $\bsig$ s.t. for any $ \xb^{(i)} \in \RR^d$  with $i\in [m]$,
	\begin{equation}\label{eq:noise}
		\frac{1}{m}\sum_{i=1}^m \EE_{\xi^{(i)}} \left[\norm{\nabla f_i(\xb^{(i)}, \xi^{(i)}) - \nabla f_i(\xb^{(i)})}^2\right] \le \bsig^2.
	\end{equation}
\end{assumption}

In this paper, we focus on the smooth and strongly-convex functions. 
That is, the function $f_i(x)$ satisfies the following assumption.
\begin{assumption}\label{ass:f}
	Each $f_i:\RR^d\to \RR$ is $\mu$-strongly convex and $L$-smooth, i.e., for any $x, y \in \RR^d$,
	\begin{align*}
		f_i(y) \ge& f_i(x) + \dotprod{\nabla f_i(x), y -x } + \frac{\mu}{2}\norm{x - y}^2,\\
		f_i(y) \le& f_i(x) + \dotprod{\nabla f_i(x), y -x } + \frac{L}{2}\norm{x - y}^2.
	\end{align*}
\end{assumption}

The agents are connected through a graph $ \mathcal{G} = \{V, E\} $ with $ V $ and $ E $ being the sets of nodes and edges. 
We assume that the graph is undirected and connected. 
$W$ is an $m\times m$ mixing matrix with $W_{i,j}$ being positive if and only if there is an edge between $i$-th and $j$-th agents. 
We also assume that $W$ satisfies the following properties.
\begin{definition}[Mixing matrix]
Matrix $W\in [0, 1]^{m\times m}$ is symmetric and doubly stochastic, that is
	$ W\mathbf{1} = \mathbf{1} $, and $\mathbf{1}^\top W = \mathbf{1}^\top$.
\end{definition}
We further suppose that the mixing matrix has the following property to achieve the information average. 
Specifically, we can represent the information exchange through matrix multiplication. 
\begin{assumption}\label{ass:mix}
	Letting $W\in\RR^{m\times m}$ be a (random) mixing matrix and parameter $\theta\in (0, 1]$, it satisfies that
\begin{equation}\label{eq:dec}
\EE_W\left[\norm{W\xb - \mathbf{1}\bbx}^2\right] \le \left(1 - \theta\right)^2\norm{\xb - \mathbf{1}\bbx}^2,\mbox{ with } \theta = 1 - \sqrt{\lambda_2\left(\EE\left[W^\top W\right]\right)}. 
\end{equation}
\end{assumption}
Assumption \ref{ass:mix} says that the mixing matrix $ W $ can achieve averaging in expectation but without any other constraint.
As a concrete example, \citet{boyd2006randomized} showed that randomized gossip matrices with time-varying topologies satisfy Assumption~\ref{ass:mix}. 

\section{Snap-Shot Decentralized Stochastic Gradient Tracking}
\label{sec:ss}

In this section, we propose an algorithm named ``Snap-Shot Decentralized Stochastic Gradient Tracking'' (\algname).
We first give the algorithm description and the intuition behind our algorithm.
Then, we provide a detailed convergence analysis of \algname. 

\begin{algorithm}[tb]
	\caption{Snap-Shot Decentralized Stochastic Gradient Tracking}
	\label{alg:alg_name}
	\small
	\begin{algorithmic}
		\STATE {\bfseries{Input}}: $x_0$, mixing matrix $W$, initial step size $\eta$.
		\STATE {\bfseries{Initialization}:} Set $\xb_0 = \mathbf{1}x_0$, $\qb_0 = \mathbf{1}x_0$, $\bs_0^{(i)} = \nabla f_i(\xb_0^{(i)},\xi_0)$, in parallel for $i \in [m]$, $\tau = 0$.
		\FOR {$t = 1,\dots, T$}
		\STATE Generate $\zeta_t$ with probability $p$. 
		\STATE Sample $\xi_t^{(i)}$ in parallel for all $m$ agents and update
		\begin{equation}\label{eq:xx_up}
			\xb_{t+1} = W \left( \xb_t - \eta_t \left( \bs_t +   \nabla F(\xb_t, \xi_t) - \nabla F(\qb_t,\xi_\tau)  \right) \right).
		\end{equation} 
		\STATE Update 
		\begin{equation}\label{eq:qq}
			\qb_{t+1} = \zeta_t\xb_t + (1 - \zeta_t) \qb_t.
		\end{equation}
		\STATE Update 
		\begin{equation} \label{eq:ss_up}
			\bs_{t+1} = W\bs_t + \zeta_t\left( \nabla F(\xb_t,\xi_t) - \nabla F(\qb_t,\xi_\tau) \right).
		\end{equation}
		\STATE Set 
		$$\tau = \begin{cases}
			t,  \quad \mbox{ if } \zeta_t = 1, \\
			\tau,\quad \mbox{otherwise}.
		\end{cases}$$  
		\ENDFOR
	\end{algorithmic}
\end{algorithm}

\subsection{Algorithm Description}

Our work extends the idea of snap-shot gradient tracking (\ssgt) proposed by \citet{song2021optimal} to the decentralized \sgd. 
The detailed algorithm description is in Algorithm~\ref{alg:alg_name}. 
%Our algorithm has the following algorithmic procedure with $ \eta_t $ being the step size
%\begin{align}
%	\xb_{t+1} =& W \Big(\xb_t - \eta_t \Big(\bs_t + \nabla F(\xb_t, \xi_t) - \nabla F(\qb_t, \xi_\tau)\Big) \Big) \label{eq:xx_up}\\
%	\qb_{t+1} =& \zeta_t \xb_t + (1 - \zeta_t) \qb_t, \quad \mbox{with} \quad \zeta_t\sim {Bernoulli}(p) \label{eq:qq}\\
%	\bs_{t+1} =& W \bs_t + \zeta_t \Big(\nabla F(\xb_t, \xi_t) - \nabla F(\qb_t, \xi_\tau)\Big) \label{eq:ss_up}\\
%	\tau =& 
%	\begin{cases}
%		t \quad &\mbox{if} \;\; \zeta_t = 1\\
%		\tau    & \mbox{else}
%	\end{cases}
%\end{align}

Following the idea of \ssgt, our algorithm introduces a variable $ \qb_t $ to record some history position of $ \xb_t $ and updates it with probability $p$.
Furthermore, instead of updating the gradient tracking variable $\bs_t$ with the aggregated stochastic gradient $\nabla F(\xb_{t+1}, \xi_{t+1})$ for each iteration in \dsgt~(refer to Eq.~\eqref{eq:ss_GT}), \algname~updates $\bs_t$ with gradient information also with probability $ p $. 
The value of $ \bbs_t $ is updated only when $ \zeta_t = 1 $.
If $\zeta_t = 1$, we need to update the $\tau$ which records the time update $\qb_t$.

Unlike \dsgt~whose $ \bs_t $ tracks the average of $ \nabla F(\xb_t, \xi_t) $ \citep{pu2021distributed}, $ \bs_t $ of our algorithm tracks the average of $ \nabla F(\qb_t, \xi_\tau) $ which is shown by the following lemma.
\begin{lemma}\label{lem:s}
Let sequence $ \{\bs_t\} $ be updated as Eq.~\eqref{eq:ss_up}. Then, for any $ 0 \le t \le T $, 	it holds that
	\begin{equation}\label{eq:s}
		\bbs_t = \frac{1}{m}\sum_{i=1}^m \nabla f_i(\qbti, \xi_\tau^{(i)}).
	\end{equation}
\end{lemma}

Since $ \bs_t $ tracks the value of $ \nabla F(\qb_t, \xi_\tau) $ instead of $ \nabla F(\xb_t, \xi_t) $, our algorithm proposes the update rule~\eqref{eq:xx_up} in contrast to using $\bs_t$ to update $ \xb_t $ directly which is used in \dsgt~(See Algorithm~\ref{alg:dsgt}). 
Such modification follows from the fact
\begin{equation*}
	\mathbf{1}^\top \Big(\bs_t + \nabla F(\xb_t, \xi_t) - \nabla F(\qb_t, \xi_\tau)\Big) 
	\stackrel{\eqref{eq:s}}{=}  \mathbf{1}^\top  \nabla F(\xb_t, \xi_t).
\end{equation*}

Note that, our algorithm is inspired by the \ssgt~proposed by \citet{song2021optimal}, and its idea originates from \texttt{SVRG} \citep{johnson2013accelerating}, \texttt{L-SVRG} \citep{kovalev2020don} and \texttt{ANITA} \citep{li2021anita}. 
However, \algname~is not an easy extension of \ssgt. 
The original \ssgt~strongly correlates to the loopless Katyusha \citep{kovalev2020don}. 
Extra variables such as $\mathbf{U}_t$ compared to our algorithm and ``negative momentum'' are important in the building of \ssgt~in \cite{song2021optimal}.
Moreover, a large part of the proof of  \ssgt~follows the framework of loopless Katyusha.
Thus, it is unknown whether the idea of \ssgt~has a broader application.
Our work tries to explore the application range of the idea \ssgt~and try to extend it to decentralized stochastic gradient descent.
%Accordingly, we propose our \algname.

%\begin{align}
%\xb_{t+1} =& W \left( \xb_t - \eta \big( \bs_t +  \nabla F(\xb_t, \xi_t) - \nabla F(\qb_t, \xi_\tau)  \big) \right)\\
%\qb_{t+1} =& \zeta_t\xb_t + (1 - \zeta_t) \qb_t \\
%\bs_{t+1} =& W\bs_t + \zeta_t\left( \nabla F(\xb_t) - \nabla F(\qb_t) \right)
%\end{align}

\subsection{Convergence Analysis}
\label{subsec:conv}

We will first study the evolution of $ \EE\left[\norm{\bPi\xb_t}^2\right] $, $ \EE\left[\norm{\bPi\bs_t}^2\right] $ and $ \EE\left[\norm{\nabla F(\qb_t) - \nabla F(\mathbf{1} x^*)}^2\right] $. 
We introduce a Lyapunov function to  describe the dynamics of consensus errors and $ \norm{\nabla F(\qb_t) - \nabla F(\mathbf{1} x^*)} $.
Let us denote by $ \cF_t $ the $ \sigma $-algebra generated by $ \{(\xi_0,\zeta_0),(\xi_1,\zeta_1),\dots,(\xi_{t-1},\zeta_{t-1})\} $ and define $ \EE\left[\cdot \mid \cF_t\right] $ as the conditional expectation given $ \cF_t $.

\begin{lemma}\label{lem:main_1}
	Suppose Assumptions~\ref{ass:noise}-\ref{ass:mix} hold.
	Let $\{\eta_t\} $ be a non-increasing sequence and satisfy $\eta_t \le \frac{\theta}{16 L}$. 
	Setting  $C_{1,t} = 4\eta_t^2 / \theta^2$, $C_{2,t} = 2\eta_t/(L\theta)$, and $p = \theta$, 
	it holds that
	\begin{equation}\label{eq:psi}
		\begin{aligned}
			\EE\left[\Psi_{t+1} \mid \cF_t\right]
			\le \left(1 - \frac{\theta}{4}\right) \cdot \Psi_t +  \frac{18m\eta_t^2}{\theta} \bsig^2 + \left( \frac{88mL\eta_t^2}{\theta} +8m\eta_t \right) \big(f(\bbx_t) - f(x^*)\big),
		\end{aligned}
	\end{equation}
	where we denote 
	\begin{equation*}
		\Psi_t \triangleq  \norm{\bPi\xb_t}^2  +  C_{1,t} \norm{\bPi\bs_t}^2 + C_{2,t}\norm{\nabla F(\qb_t) - \nabla F(\mathbf{1}x^*)}^2.
	\end{equation*}
\end{lemma}

Lemma~\ref{lem:main_1} shows that $ \Psi_t $ will converge to zero under the condition that the step size $ \eta_t $ will decrease to zero and $ f(\bbx_t) - f(x^*) $ is non-increasing.
This implies that $ \norm{\bPi\xb_t} $ will converge to zero, that is, the distance $ \norm{\xbti - \bbx_t}^2  $ will vanish as $ t $ goes.
Next, we are going to upper bound  the distance $ \norm{\bbx_t - x^*} $.

\begin{lemma}\label{lem:main_2}
Suppose Assumptions~\ref{ass:noise}-\ref{ass:mix} hold. Then we have the following inequality:
\begin{equation}\label{eq:xx_a}
\begin{aligned}
\EE\left[ \norm{\bbx_{t+1} - x^*}^2 \mid \cF_t\right]
\le& \left(1 - \frac{\mu\eta_t}{2}\right) \norm{\bbx_t - x^*}^2 - 2\eta_t \left(1 - 2\eta_t L \right)\big(f(\bbx_t) - f(x^*)\big) \\
&+ \eta_t^2\cdot \frac{\bsig^2}{m}
+\frac{2L\eta_t \left( 1 + 2\eta_t L \right) }{m} \norm{\bPi \xb_t}^2.
\end{aligned}
\end{equation}
\end{lemma}

Lemma~\ref{lem:main_1} and~\ref{lem:main_2} show that the dynamics of $ \Psi_t $ and $ \norm{\bbx_t - x^*}^2 $ correlate to each other.
Based on above two lemmas, we obtain the following convergence properties.

\begin{lemma}\label{lem:main_3}
Suppose Assumptions~\ref{ass:noise}-\ref{ass:mix} hold.
Let $\{\eta_t\} $ be a non-increasing sequence and satisfy $\eta_t \le \frac{\theta}{2^6\cdot 3 \cdot L}$. 
It holds that
\begin{equation}\label{eq:xx_dec}
\begin{aligned}
&\EE\left[\norm{\bbx_{t+1} - x^*}^2 + \frac{24L\eta_{t+1}}{m\theta} \Psi_{t+1} \mid \cF_t\right]\\
\le&
\exp\left(- \frac{\mu\eta_t}{2}\right) \left(\norm{\bbx_t - x^*}^2 +  \frac{24L\eta_t}{m\theta} \Psi_t \right) 
- \frac{7\eta_t}{8} \big(f(\bbx_t) - f(x^*)\big) 
+ \frac{2^4 \cdot 3^3 \cdot L\eta_t^3}{\theta^2}\bsig^2 + \eta_t^2\cdot \frac{\bsig^2}{m}.	
\end{aligned}
\end{equation}
\end{lemma}

Based on Lemma~\ref{lem:main_3}, we can derive the desired convergence properties shown in the following theorem.  

\begin{theorem}\label{thm:main}
Suppose Assumptions~\ref{ass:noise}-\ref{ass:mix} hold. Sequences $ \{\xb_t\} $, $ \{\qb_t\} $, and $ \{\bs_t\} $ are generated by Algorithm~\ref{alg:alg_name}. 
Then Algorithm~\ref{alg:alg_name} has the following convergence properties:
\begin{itemize}
	\item If $ \bsig^2 = 0 $ and step size $ \eta_t =  \frac{\theta}{2^6\cdot 3 \cdot L}$, it holds that
	\begin{equation}\label{eq:c1}
		\EE\left[\norm{\bbx_T - x^*}^2 + \frac{1}{8m} \Psi_T\right]
		\le \exp\left( -\frac{1}{2^7\cdot 3} \cdot \frac{\theta\mu}{L} \cdot  T \right) \cdot \left(\norm{\bbx_0 - x^*}^2 + \frac{1}{8m} \Psi_0\right).
	\end{equation}
\item If $ \bsig^2 >0 $, set the step size sequence $ \{\eta_t\} $ and weight sequence $ \{\omega_t\} $ as follows:
\begin{align*}
	\eta_t = \frac{6\beta}{L + \beta \mu t},\quad \mbox{ and } \quad \omega_t = \frac{\eta_t}{\eta_0} \exp\left(\frac{\mu}{2} \sum_{i=0}^t \eta_i \right),\quad \mbox{with} \quad\beta = \frac{\theta}{2^7\cdot 3^2}.
\end{align*}
Letting $S_T = \sum_{t=0}^T \omega_t$, then it holds that
\begin{equation}\label{eq:c2}
\begin{aligned}
	&\frac{1}{S_T} \sum_{t=0}^T \omega_t \big(\EE  [f(\bbx_t)] - f(x^*)\big) \\
	\le&
	\frac{24 \cdot (L + \beta\mu (T+1))^2}{ \beta^2\mu^3 T^3} \cdot \frac{\bsig^2}{m}
	+ \frac{2^8\cdot 3^5 \cdot L}{\mu^2\theta^2} \cdot \frac{\bsig^2}{T^2} + \frac{2^{24}\cdot 3^6\cdot  L^3}{\theta^3 \mu^2} \cdot \frac{1}{T^3} \cdot \left(\norm{\bbx_0 - x^*}^2 + \frac{\Psi_0}{8m}\right).
\end{aligned}
\end{equation}
\end{itemize}
\end{theorem}
\begin{proof}
	For the case $ \bsig^2 = 0 $, Eq.~\eqref{eq:xx_dec} reduces to 
	\begin{align*}
	&\EE\left[\norm{\bbx_{t+1} - x^*}^2 + \frac{24L\eta_{t+1}}{m\theta} \Psi_{t+1} \mid \cF_t\right]\\
	\le&  
	\exp\left(- \frac{\mu\eta_t}{2}\right) \left(\norm{\bbx_t - x^*}^2 +  \frac{24L\eta_t}{m\theta} \Psi_t \right) 
	- \frac{7\eta_t}{8} \big(f(\bbx_t) - f(x^*)\big) \\
	\le&
	\exp\left(- \frac{\mu\eta_t}{2}\right) \left(\norm{\bbx_t - x^*}^2 +  \frac{24L\eta_t}{m\theta} \Psi_t \right).
	\end{align*}
	Using above equation recursively and replacing $ \eta_t =  \frac{\theta}{2^{6}\cdot 3\cdot L}$, we can obtain the first result.

The result for the case $ \bsig^2 >0 $ follows from Lemma~\ref{lem:av} with $e_t = f(\bbx_t) - f(x^*)$, $r_t = \norm{\bbx_t - x^*}^2 +  \frac{24L\eta_t}{m\theta} \Psi_t$, $ A = \frac{7}{8} $, $ B = \frac{\bsig^2}{m}$, and $ C =\frac{2^4 \cdot 3^3 \cdot L}{\theta^2}\bsig^2 $.
Specifically, we have
\begin{align*}
	&\frac{1}{S_T} \sum_{t=0}^T \omega_t \big(\EE  [f(\bbx_t)] - f(x^*)\big) \\
	\stackrel{\eqref{eq:SS}}{\le}& 
	\frac{2^{24}\cdot 3^6\cdot  L^3}{\theta^3 \mu^2} \cdot \frac{1}{T^3} \cdot \left(\norm{\bbx_0 - x^*}^2 + \frac{24L\eta_0}{m\theta} \Psi_0\right) + \frac{18\cdot 8 \cdot L^3 (L + \beta\mu (T+1))^2}{7 \cdot \beta^2(L - \beta \mu)^3 \mu^3 T^3} \cdot \frac{\bsig^2}{m}  \\
	&+ \frac{2^9\cdot 3^6\cdot L^4}{7 \mu^2 \theta^2(L - \beta \mu)^3} \cdot \frac{\bsig^2}{T^2}\\
	\le&
	\frac{2^{24}\cdot 3^6\cdot  L^3}{\theta^3 \mu^2} \cdot \frac{1}{T^3} \cdot \left(\norm{\bbx_0 - x^*}^2 + \frac{\Psi_0}{8m}\right) 
	+ \frac{18\cdot 8 \cdot L^3 (L + \beta\mu (T+1))^2}{7 \cdot (6L^3/7) \cdot \beta^2\mu^3 T^3} \cdot \frac{\bsig^2}{m}  + \frac{2^9\cdot 3^6\cdot L^4}{7 \mu^2 \theta^2 \cdot (6L^3/7)} \cdot \frac{\bsig^2}{T^2}\\
	=&
	\frac{2^{24}\cdot 3^6\cdot  L^3}{\theta^3 \mu^2} \cdot \frac{1}{T^3} \cdot \left(\norm{\bbx_0 - x^*}^2 + \frac{\Psi_0}{8m}\right)
	+\frac{24 \cdot (L + \beta\mu (T+1))^2}{ \beta^2\mu^3 T^3} \cdot \frac{\bsig^2}{m}
	+ \frac{2^8\cdot 3^5 \cdot L}{\mu^2\theta^2} \cdot \frac{\bsig^2}{T^2},
\end{align*}
where the second inequality is because of $0<\theta < 1$ and $\mu\le L$ and 
\begin{align*}
	L - \beta \mu = L - \frac{\theta}{2^7\cdot 3^2}\mu \ge  L - \frac{L}{2^7\cdot 3^2} \ge L \left(\frac{6}{7}\right)^{1/3}.
\end{align*} 

%Using the definition of $\Psi_0$ and the value of $\eta_0$, we can obtain that
%\begin{align*}
%	\Psi_0 
%	=& 
%	\norm{\bPi\xb_0}^2 +  C_{1,0} \norm{\bPi\bs_t}^2 + C_{2,0}\norm{\nabla F(\qb_t) - \nabla F(\mathbf{1}x^*)}^2 \\
%	=& 
%	C_{1,0} \sum_{i=1}^m\norm{\nabla f_i(\bbx_0, \xi_0^{(i)}) - \frac{1}{m}\sum_{j=0}^m \nabla f_j(\bbx_0, \xi_0^{(j)})}^2 + C_{2,0} \sum_{i=1}^m \norm{\nabla f_i(\bbx_0) - \nabla f_i(x^*)}^2\\
%	=& 
%	\frac{1}{2^8\cdot 3^2\cdot L^2}\sum_{i=1}^m\norm{\nabla f_i(\bbx_0, \xi_0^{(i)}) - \frac{1}{m}\sum_{j=0}^m \nabla f_j(\bbx_0, \xi_0^{(j)})}^2\\
%	&+\frac{1}{2^5\cdot 3\cdot L^2} \sum_{i=1}^m \norm{\nabla f_i(\bbx_0) - \nabla f_i(x^*)}^2,
%\end{align*}
%where the second equality is because of initialization of Algorithm~\ref{alg:alg_name}.
%
%Therefore, 
%\begin{align*}
%&\frac{1}{S_T} \sum_{t=0}^T \omega_t \big(\EE  [f(\bbx_t)] - f(x^*)\big)\\
%\le& 
%\frac{2^7\cdot 3\cdot  L^3}{\theta^3 \mu^2} \cdot \frac{1}{T^3} \cdot \norm{\bbx_0 - x^*}^2 + \frac{L}{48m\mu^2 \theta^3} \sum_{i=1}^m\norm{\nabla f_i(\bbx_0, \xi_0^{(i)}) - \frac{1}{m}\sum_{j=0}^m \nabla f_j(\bbx_0, \xi_0^{(j)})}^2
%\end{align*}
\end{proof}

Based on Theorem~\ref{thm:main}, we can  directly obtain  the following iteration complexities.
\begin{corollary}\label{cor:complexity}
Suppose Assumptions~\ref{ass:noise}-\ref{ass:mix} hold and the mixing matrix $ W $ is fixed.  Parameters of Algorithm~\ref{alg:alg_name} are set as Theorem~\ref{thm:main}, then Algorithm~\ref{alg:alg_name} has the following iteration complexity:
\begin{itemize}
	\item If $ \bsig^2 = 0 $, to achieve $ \varepsilon $-suboptimality, the iteration complexity of Algorithm~\ref{alg:alg_name} is 
	\begin{equation}
		T = \cO\left(\frac{L}{\mu (1-\lambda_2(W))} \log \frac{1}{\varepsilon}\right)
	\end{equation}
\item If $ \bsig^2 = 0 $, to achieve $ \varepsilon $-suboptimality, the iteration complexity of Algorithm~\ref{alg:alg_name} is
\begin{equation}\label{eq:T_2}
	T = \cO\left( \frac{\bsig^2}{\mu m \varepsilon} + \frac{\sqrt{L}\bsig}{ (1 - \lambda_2(W))\mu \sqrt{\varepsilon}   } + \frac{L}{(1 - \lambda_2(W)) \mu \varepsilon^{1/3}} \right)
\end{equation}
\end{itemize}
\end{corollary}
\begin{proof}
If $ \bsig^2 = 0 $, to achieve $ \varepsilon $-suboptimality, by Eq.~\eqref{eq:c1}, it requires that
\begin{align*}
\exp\left( -\frac{1}{2^7\cdot 3} \cdot \frac{\theta\mu}{L} \cdot  T \right) \cdot \left(\norm{\bbx_0 - x^*}^2 + \frac{1}{8m} \Psi_0\right) \le \varepsilon,
\end{align*}
which leads to 
\begin{align*}
	T = \frac{2^7\cdot 3 \cdot L}{\mu \theta} \log \frac{\norm{\bbx_0 - x^*}^2 + \frac{1}{8m} \Psi_0}{\varepsilon} = \cO\left(\frac{L}{\mu (1 - \lambda_2(W))} \log \frac{1}{\varepsilon}\right),
\end{align*}
where the last equality is because of the fact that $ \theta = 1 - \lambda_2(W) $ when $ W $ is a fixed mixing matrix.

If $ \bsig^2 > 0  $, supposing $ T $ is sufficient large that $ \beta \mu (T+1) $ dominates $ L $, in this case, Eq.~\eqref{eq:c2} reduces to
\begin{align*}
&\frac{1}{S_T} \sum_{t=0}^T \omega_t \big(\EE  [f(\bbx_t)] - f(x^*)\big) \\
=&
\cO\left( \frac{1}{T} \cdot \frac{\bsig^2}{m\mu} + \frac{L}{\mu^2 \theta^2} \cdot \frac{\bsig^2}{T^2} + \frac{L^3}{\theta^3 \mu ^3} \cdot  \frac{1}{T^3} \cdot \left(\norm{\bbx_0 - x^*}^2 + \frac{\Psi_0}{8m}\right)\right).
\end{align*}
Thus, to achieve $ \varepsilon $-suboptimality, the iteration complexity is 
\begin{align*}
	T = \cO\left( \frac{\bsig^2}{\mu m \varepsilon} + \frac{\sqrt{L}\bsig}{ \mu \theta \sqrt{\varepsilon}  } + \frac{L}{\theta \mu \varepsilon^{1/3}} \right).
\end{align*}
Replacing  $ \theta = 1 - \lambda_2(W) $ to above equation concludes the proof.
\end{proof}

\begin{remark}

\begin{itemize}
    \item When $ \bsig^2 = 0 $, our algorithm achieves a linear convergence rate and its iteration complexity is $ \cO\left(\frac{L}{\mu (1 - \lambda_2(W))} \log \frac{1}{\varepsilon}\right) $.
In contrast, the iteration complexity of \dsgt~is $ \cO\left(\frac{L}{\mu(1 - \lambda_2(W))^2} \log\frac{1}{\varepsilon}\right) $ \citep{pu2021distributed,qu2017harnessing}. Thus, our \algname~has better performance than \dsgt~theoretically.

\item When $ \bsig^2 > 0 $, we can observe that the iteration complexity of \algname~depends on $ (1 - \lambda_2(W))^{-1} $ while \dsgt~depends on $ (1 - \lambda_2(W))^{-3/2} $ (see Eq.~\eqref{eq:com1}). 
Thus, \algname~also outperforms \dsgt~when $\bsig^2 >0$.
\end{itemize}

\end{remark}

\section{Acceleration with Loopless Chebyshev Acceleration}
\label{sec:ass}

In this section, we try to further improve \algname~and combine it with  the loopless Chebyshev acceleration proposed by \citet{song2021optimal}.
Because the loopless Chebyshev acceleration only works for the static networks, we assume that for all iterations, it shares the same mixing matrix $ W $ in this section. 
For the static networks, Corollary~\ref{cor:complexity} shows that the iteration complexity of \algname~depends on $ \theta^{-1} = (1 - \lambda_2(W))^{-1} $.
In this section, we propose \algnamea~which achieves an iteration complexity depending on $ \Big(1 - \lambda_2(W)\Big)^{-1/2} $ instead of $ \Big(1 - \lambda_2(W)\Big)^{-1} $. 

\subsection{Algorithm Description}

Before introducing \algnamea, we make an additional assumption on the mixing matrix and define some new necessary notations.
\begin{assumption}\label{ass:mix2}
	The mixing matrix $ W\in\RR^{m\times m} $ is positive semi-definite.
\end{assumption}
The above assumption can be easily satisfied since we can choose $ \frac{I + W}{2} $ as the mixing matrix that is positive semi-definite for any mixing matrix $ W $.

Now, we introduce $ 2m\times 2m $ augmented matrices  $ \tW $ and $ \tPi $  for the  mixing matrix $ W $ and projection matrix $ \bPi $ defined as follows:
\begin{equation}
	\tW = \begin{bmatrix}
		(1 + \gamma) W & -\gamma W \\
		I_m & \mathbf{0}
	\end{bmatrix}, \quad \mbox{and} \quad \tPi = \begin{bmatrix}
	\bPi & \mathbf{0}\\
	\mathbf{0} & \bPi
\end{bmatrix}.
\end{equation}
Accordingly, we define the augmented decision variable $ \txb \in\RR^{2m\times d}$ and gradient-tracking variable $ \tsb \in\RR^{2m\times d} $.
Furthermore, we denote that $ \xb_t := \txb_t^{(1:m)} $, that is, $ \xb_t $ takes the value of the first $ m $ rows of $ \txb_t $.
Given these notations, we describe \algnamea~in Algorithm~\ref{alg:alg_namea}.

\begin{algorithm}[tb]
	\caption{Snap-Shot Decentralized Stochastic Gradient Tracking with Loopless Chebyshev Acceleration (\algnamea)}
	\label{alg:alg_namea}
	\small
	\begin{algorithmic}
		\STATE {\bfseries{Input}}: $x_0$, mixing matrix $W$, initial step size $\eta$.
		\STATE {\bfseries{Initialization}:} Set $\txb_0 = [\mathbf{1}x_0; \mathbf{1}x_0]$, $\qb_0 = \mathbf{1}x_0$, $\bs_0^{(i)} = \nabla f_i(\xb_0^{(i)},\xi_0)$, in parallel for $i \in [m]$, $ \tsb_{0} = [\bs_0; \bs_0] $, and $\tau = 0$.
		\FOR {$t = 1,\dots, T$}
		\STATE Generate $\zeta_t$ with probability $p$. 
		\STATE Sample $\xi_t^{(i)}$ in parallel for all $m$ agents and update
		\begin{equation}\label{eq:xx_up_1}
			\txb_{t+1} = \tW\left( \txb_t - \eta_t \Big( \tsb_t +\left[\nabla F(\xb_t, \xi_t);\;\nabla F(\xb_t, \xi_t)  \right] - 
			\left[\nabla F(\qb_t, \xi_\tau);\; \nabla F(\qb_t, \xi_\tau)  \right]
			\Big)\right).
		\end{equation} 
		\STATE Update 
		\begin{equation}\label{eq:qq_1}
			\qb_{t+1} = \zeta_t\xb_t + (1 - \zeta_t) \qb_t.
		\end{equation}
		\STATE Update 
		\begin{equation} \label{eq:ss_up_1}
			\tsb_{t+1} = \tW\tsb_t + \zeta_t \bigl( \left[\nabla F(\xb_t, \xi_t);\;\nabla F(\xb_t, \xi_t)  \right] - 
			\left[\nabla F(\qb_t, \xi_\tau);\; \nabla F(\qb_t, \xi_\tau)  \right] \bigr).
		\end{equation}
		\STATE Set 
		$$\tau = \begin{cases}
			t,  \quad \mbox{ if } \zeta_t = 1, \\
			\tau,\quad \mbox{otherwise}.
		\end{cases}$$  
		\ENDFOR
	\end{algorithmic}
\end{algorithm}

We can observe that Algorithm~\ref{alg:alg_namea} shares almost the same algorithmic structure to the one of Algorithm~\ref{alg:alg_name}.
The advantage of \algnamea~mainly depends on the following property.
\begin{lemma}[Lemma 11 of \citet{song2021optimal}]
	Under Assumption~\ref{ass:mix2}, for any $ \xb\in\RR^{m\times d} $ and $ t>0 $, it holds that
	\begin{align}
		\norm{\tPi \tW^t \tPi [\xb;\;\xb]}^2 \le \alpha \left(1 - \tilde{\theta}\right)^{2t}\norm{\bPi \xb}^2, \quad\mbox{with}\quad \alpha \le 14 \quad\mbox{and}\quad \tilde{\theta} = \cO\left( \sqrt{1 -\lambda_2(W)}\right). \label{eq:tW_dec}
	\end{align}
\end{lemma}
The above property is also used in the analysis of the heavy ball method and shows that the heavy ball method can achieve a faster convergence rate than the gradient descent \citep{recht2010cs726}.
%\begin{align}
%\txb_{t+1} =& \tW\left( \txb_t - \eta \left( \tsb_t + \begin{bmatrix}
%\nabla F(\xb_t, \xi_t)\\
%\nabla F(\xb_t, \xi_t)
%\end{bmatrix} - 
%\begin{bmatrix}
%	\nabla F(\qb_t, \xi_\tau)\\
%	\nabla F(\qb_t, \xi_\tau)
%\end{bmatrix}  \right)\right)\\
%\qb_{t+1} =& \zeta_t \xb_t + (1 - \zeta_t)\qb_t\\
%\tsb_{t+1} =& \tW\tsb_t + \zeta_t \bigl( \nabla F(\xb_t, \xi_t)_{\#} - \nabla F(\qb_t, \xi_\tau)_{\#} \bigr)
%\end{align}

\subsection{Convergence Analysis}

First, we will show that the first and last $ m $ rows of $ \txb_t $ share the same mean. 
This property also holds for $ \tsb_t $. 
\begin{lemma}\label{lem:equal}
	Letting sequences $ \{\txb_t\} $ and $ \{\tsb_t\} $ are generated by Algorithm~\ref{alg:alg_namea}, then it holds that
	\begin{align}
		\sum_{i=1}^{m} \txb_t^{(i)} = \sum_{i=m+1}^{2m} \txb_t^{(i)}, \quad  \sum_{i=1}^{m} \tsb_t^{(i)} = \sum_{i=m+1}^{2m} \tsb_t^{(i)}, \label{eq:equal}
	\end{align}
	and
	\begin{equation} \label{eq:sq}
		\frac{1}{2m}\sum_{i=1}^{2m} \tsb_t^{{i}} = \frac{1}{m}\sum_{i=1}^m \nabla f_i(\qbti, \xi_\tau^{(i)}).
	\end{equation}
\end{lemma}

Above lemma shows that the means of $ \txb_t $ and $ \tsb_t $ equal to $ \bbx_t $ and $ \bbs_t $, respectively.
Thus, Lemma~\ref{lem:main_2} still holds for \algnamea. 
Next, we will focus on analyzing the convergence properties of consensus errors which are different from the ones of \algname.

\begin{lemma}\label{lem:tsx}
Letting sequences $ \{\txb_t\} $ and $ \{\tsb_t\} $ are generated by Algorithm~\ref{alg:alg_namea}, it holds that
\begin{align}
\EE\left[ \norm{\tPi\txb_t}^2 \right] \le \cC_{x,t} \quad \mbox{and}\quad \EE\left[\norm{\tPi\tsb_t}^2\right] \le \cC_{s,t}
 \label{eq:tsx} 
\end{align}
with 
\begin{align}
\cC_{s,t+1} \le& \left(1-\frac{\tilde{\theta}}{2}\right)\cC_{s,t} + 4\alpha p \left(\norm{\nabla F(\xb_t) - \nabla F(\qb_t)}^2 + 2m\bsig^2\right) \label{eq:ts} \\
\cC_{x,t+1} \le&  \left(1 - \frac{\tilde{\theta}}{2}\right) \cC_{x,t} +  \frac{3\alpha}{\tilde{\theta}}  \eta_t^2 \left( \cC_{s,t} + 2 \norm{\nabla F(\xb_t) - \nabla F(\qb_t)}^2  \right) + 2\alpha m \bsig^2  \eta_t^2 \label{eq:tx}
\end{align}
and 
\begin{align}
\cC_{s,0} = 2\alpha \norm{\Pi \bs_0}^2, \quad \cC_{x,0} = \alpha \norm{\tPi\xb_0}^2.
\end{align}
\end{lemma}

Based on the above lemma about the consensus error terms, we can obtain the following lemma similar to Lemma~\ref{lem:main_1}.
\begin{lemma}\label{lem:psi_1}
	Suppose Assumptions~\ref{ass:noise}-\ref{ass:f} and Assumption~\ref{ass:mix2} hold.
	Let $\{\eta_t\} $ be a non-increasing sequence and satisfy $ \eta_t \le \frac{\tilde{\theta}}{2^4\cdot 3^3\cdot L} $. 
	Setting $ C_{1, t} = \frac{12\alpha\eta_t^2}{\tilde{\theta}^2} $, $ C_{2,t} = \frac{16(1+8\alpha)\eta_t^2}{\tilde{\theta}^2} $, $ p = \tilde{\theta} $, then we can obtain that
	\begin{align} \label{eq:psi_1}
		\widetilde{\Psi}_{t+1} 
		\le 
		\left(1 - \tilde{\theta}\right)\widetilde{\Psi}_t + \frac{2^{12}\cdot 3^2\cdot m \eta_t^2}{\tilde{\theta}} \cdot \bsig^2 + \frac{2^{12}\cdot 3^2\cdot m L \eta_t^2}{\tilde{\theta}} \cdot \Big(f(\bbx_t) - f(x^*)\Big),
	\end{align}
	where we define
	\begin{align}
		\widetilde{\Psi}_t := \cC_{x,t} +  C_{1,t} \cdot \cC_{s, t} + C_{2,t} \norm{\nabla F(\qb_t) - \nabla F(\mathbf{1}x^*)}^2.
	\end{align}
\end{lemma}

Combining Lemma~\ref{lem:psi_1} and Lemma~\ref{lem:main_2}, we can obtain a lemma similar to Lemma~\ref{lem:main_3}.
\begin{lemma}\label{lem:main_tld}
	Suppose Assumptions~\ref{ass:noise}-\ref{ass:f} and Assumption~\ref{ass:mix2} hold.
	Let $\{\eta_t\} $ be a non-increasing sequence and satisfy $ \eta_t \le \frac{\tilde{\theta}}{2^8\cdot 3\cdot L} $. Then, it holds that
\begin{equation}\label{eq:ttt}
	\begin{aligned}
		&\EE\left[ \norm{\bbx_{t+1} - x^*}^2 + \frac{48L\eta_{t+1}}{m\tilde{\theta}} \widetilde{\Psi}_{t+1} \right]
		\le 
		\exp\left(-\frac{\mu\eta_t}{2}\right)\left(\norm{\bbx_t - x^*}^2 + \frac{48L\eta_t}{m\tilde{\theta}} \cdot \widetilde{\Psi}_t\right) \\
		&- \frac{7\eta_t}{8} \Big(f(\bbx_t) - f(x^*)\Big) + \eta_t^2\cdot \frac{\bsig^2}{m} + \frac{2^{16}\cdot 3^2 \cdot  L\eta_t^3}{\tilde{\theta}^2} \cdot \bsig^2.
	\end{aligned}
\end{equation}
\end{lemma}

\begin{theorem}\label{thm:main_2}
	Suppose Assumptions~\ref{ass:noise}-\ref{ass:f} and Assumption~\ref{ass:mix2} hold. Sequences $ \{\xb_t\} $, $ \{\qb_t\} $, and $ \{\tsb_t\} $ are generated by Algorithm~\ref{alg:alg_namea}. 
	Then Algorithm~\ref{alg:alg_namea} has the following convergence properties:
	\begin{itemize}
		\item If $ \bsig^2 = 0 $ and step size $ \eta_t =  \frac{\theta}{2^8\cdot 3 \cdot L}$, it holds that
		\begin{equation}\label{eq:cc1}
			\EE\left[\norm{\bbx_T - x^*}^2 + \frac{1}{8m} \Psi_T\right]
			\le \exp\left( -\frac{1}{2^7\cdot 3} \cdot \frac{\theta\mu}{L} \cdot  T \right) \cdot \left(\norm{\bbx_0 - x^*}^2 + \frac{1}{8m} \Psi_0\right).
		\end{equation}
		\item If $ \bsig^2 >0 $, set the step size sequence $ \{\eta_t\} $ and weight sequence $ \{\omega_t\} $ as follows:
		\begin{align*}
			\eta_t = \frac{6\tbeta}{L + \tbeta \mu t},\quad \mbox{ and } \quad \omega_t = \frac{\eta_t}{\eta_0} \exp\left(\frac{\mu}{2} \sum_{i=0}^t \eta_i \right),\quad \mbox{with} \quad\tbeta = \frac{\tilde{\theta}}{2^9\cdot 3^2}.
		\end{align*}
		Letting $S_T = \sum_{t=0}^T \omega_t$, then it holds that
		\begin{equation}\label{eq:cc2}
			\begin{aligned}
				&\frac{1}{S_T} \sum_{t=0}^T \omega_t \big(\EE  [f(\bbx_t)] - f(x^*)\big) \\
				\le&
				\frac{24 \cdot (L + \tbeta\mu (T+1))^2}{ \tbeta^2\mu^3 T^3} \cdot \frac{\bsig^2}{m}
				+ \frac{2^{20}\cdot 3^4 \cdot L}{\tilde{\theta}^2\mu^2} \cdot \frac{\bsig^2}{T^2} + \frac{2^{28}\cdot 3^5\cdot  L^3}{\tilde{\theta}^3 \mu^2} \cdot \frac{1}{T^3} \cdot \left(\norm{\bbx_0 - x^*}^2 + \frac{\widetilde{\Psi}_0}{16m}\right).
			\end{aligned}
		\end{equation}
	\end{itemize}
\end{theorem}
\begin{proof}
	For the case $ \bsig^2 = 0 $, Eq.~\eqref{eq:ttt} reduces to 
	\begin{align*}
		&\EE\left[\norm{\bbx_{t+1} - x^*}^2 + \frac{48L\eta_{t+1}}{m\tilde{\theta}} \widetilde{\Psi}_{t+1}\right]\\
		\le&  
		\exp\left(- \frac{\mu\eta_t}{2}\right) \left(\norm{\bbx_t - x^*}^2 +  \frac{48L\eta_t}{m\tilde{\theta}} \widetilde{\Psi}_t \right) 
		- \frac{7\eta_t}{8} \big(f(\bbx_t) - f(x^*)\big) \\
		\le&
		\exp\left(- \frac{\mu\eta_t}{2}\right) \left(\norm{\bbx_t - x^*}^2 +  \frac{24L\eta_t}{m\tilde{\theta}} \widetilde{\Psi}_t \right).
	\end{align*}
	Using above equation recursively and replacing $ \eta_t =  \frac{\theta}{2^{8}\cdot 3\cdot L}$, we can obtain the first result.

	The result for the case $ \bsig^2 >0 $ follows from Lemma~\ref{lem:av} with $e_t = f(\bbx_t) - f(x^*)$, $r_t = \norm{\bbx_t - x^*}^2 +  \frac{48L\eta_t}{m\tilde{\theta}} \Psi_t$, $ A = \frac{7}{8} $, $ B = \frac{\bsig^2}{m}$, and $ C =\frac{2^{16} \cdot 3^3 \cdot L}{\tilde{\theta}^2}\bsig^2 $.
	Specifically, we have
	\begin{align*}
		&\frac{1}{S_T} \sum_{t=0}^T \omega_t \big(\EE  [f(\bbx_t)] - f(x^*)\big) \\
		\stackrel{\eqref{eq:SS}}{\le}& 
		\frac{2^{28}\cdot 3^5\cdot  L^3}{\tilde{\theta}^3 \mu^2} \cdot \frac{1}{T^3} \cdot \left(\norm{\bbx_0 - x^*}^2 + \frac{48L\eta_0}{m\tilde{\theta}} \widetilde{\Psi}_0\right) + \frac{18\cdot 8 \cdot L^3 (L + \tbeta\mu (T+1))^2}{7 \cdot \tbeta^2(L - \tbeta \mu)^3 \mu^3 T^3} \cdot \frac{\bsig^2}{m}  \\
		&+ \frac{2^{21}\cdot 3^5\cdot L^4}{7 \mu^2 \theta^2(L - \tbeta \mu)^3} \cdot \frac{\bsig^2}{T^2}\\
		\le&
		\frac{2^{28}\cdot 3^5\cdot  L^3}{\tilde{\theta}^3 \mu^2} \cdot \frac{1}{T^3} \cdot \left(\norm{\bbx_0 - x^*}^2 + \frac{\widetilde{\Psi}_0}{16m}\right) 
		+ \frac{18\cdot 8 \cdot L^3 (L + \tbeta\mu (T+1))^2}{7 \cdot (6L^3/7) \cdot \tbeta^2\mu^3 T^3} \cdot \frac{\bsig^2}{m}  + \frac{2^{21}\cdot 3^5\cdot L^4}{7 \mu^2 \theta^2 \cdot (6L^3/7)} \cdot \frac{\bsig^2}{T^2}\\
		=&
		\frac{2^{28}\cdot 3^5\cdot  L^3}{\tilde{\theta}^3 \mu^2} \cdot \frac{1}{T^3} \cdot \left(\norm{\bbx_0 - x^*}^2 + \frac{\widetilde{\Psi}_0}{16m}\right)
		+\frac{24 \cdot (L + \tbeta\mu (T+1))^2}{ \tbeta^2\mu^3 T^3} \cdot \frac{\bsig^2}{m}
		+ \frac{2^{20}\cdot 3^4 \cdot L}{\mu^2\tilde{\theta}^2} \cdot \frac{\bsig^2}{T^2},
	\end{align*}
	where the second inequality is because of $0<\theta < 1$ and $\mu\le L$ and 
	\begin{align*}
		L - \tbeta \mu = L - \frac{\theta}{2^7\cdot 3^2}\mu \ge  L - \frac{L}{2^7\cdot 3^2} \ge L \left(\frac{6}{7}\right)^{1/3}.
	\end{align*} 
\end{proof}

\begin{corollary}\label{cor:complexity_2}
	Suppose Assumptions~\ref{ass:noise}-\ref{ass:f} and Assumption~\ref{ass:mix2} hold. Parameters of Algorithm~\ref{alg:alg_namea} are set as Theorem~\ref{thm:main_2}, then Algorithm~\ref{alg:alg_namea} has the following iteration complexity:
	\begin{itemize}
		\item If $ \bsig^2 = 0 $, to achieve $ \varepsilon $-suboptimality, the iteration complexity of Algorithm~\ref{alg:alg_name} is 
		\begin{equation}
			T = \cO\left(\frac{L}{\mu \sqrt{1 - \lambda_2(W)}} \log \frac{1}{\varepsilon}\right)
		\end{equation}
		\item If $ \bsig^2 > 0 $, to achieve $ \varepsilon $-suboptimality, the iteration complexity of Algorithm~\ref{alg:alg_name} is
		\begin{equation}\label{eq:T_3}
			T = \cO\left( \frac{\bsig^2}{\mu m \varepsilon} + \frac{\sqrt{L}\bsig}{ \mu \sqrt{1-\lambda_2(W)} \sqrt{\varepsilon}  } + \frac{L}{\sqrt{1-\lambda_2(W)} \mu \varepsilon^{1/3}} \right)
		\end{equation}
	\end{itemize}
\end{corollary}

\begin{remark}
	Eq.~\eqref{eq:T_3} shows that the loopless Chebyshev acceleration can effectively reduce the iteration complexity.
	Comparing Eq.~\eqref{eq:T_3} with Eq.~\eqref{eq:com2}, we can conclude that \algnamea~can achieve good performance comparable to \dsgt~for all cases (no $C_W$).
	In contrast,  \dsgt~can only achieve good performance on the limited kinds of communication networks.
\end{remark}

\section{Conclusion}

In this paper, we explore the application range of the idea of \ssgt~and extend it to design novel decentralized  \texttt{SGD} methods.
We propose two novel algorithms named \algname~and \algnamea~based on the idea of \ssgt.
These two algorithms have similar algorithmic structure to \dsgt~and they both take single loop communication strategy, which is the same as \dsgt. 
Both \algname~and \algnamea~can achieve better convergence rate than \dsgt~in general cases. 
Especially, \algnamea~achieves the best convergence rate among decentralized \texttt{SGD} methods  with  single loop communication strategy.

\bibliography{ref}
\bibliographystyle{plainnat}

\appendix
\begin{appendix}
	\onecolumn
	\begin{center}
		{\huge {Supplementary Material to "Snap-Shot Decentralized Stochastic Gradient Tracking Methods"}}
	\end{center}

\section{Useful Lemmas}

In this section, we will introduce several useful lemmas that will be used in our proofs.
These lemmas are easy to check or prove. 
Thus, we omit the detailed proofs of these lemmas. 
\begin{lemma}
Let $g(x)$ be a monotonically increasing function in the range $[t_0,T]$, then it holds
that
\begin{equation} \label{eq:int_1}
\int_{t_0}^{T} g(x)\; dx \le	\sum_{k=t_0}^T g(k) \le \int_{t_0}^{T+1} g(x)\; dx.
\end{equation}
If f(x) is monotonically decreasing in the range  $[t_0,T]$, then it holds that
\begin{equation}\label{eq:int_2}
\int_{t_0}^{T+1} g(x)\; dx \le	\sum_{k=t_0}^T g(k) \le \int_{t_0-1}^{T} g(x)\; dx.
\end{equation}
\end{lemma}

\begin{lemma}
	If $a_i$'s are independent random variables with expectation $\EE[a_i] = 0$, then it holds that
	\begin{equation}\label{eq:ee}
		\EE\left[\norm{\frac{1}{m}\sum_{i=1}^m a_i }^2\right] = \frac{1}{m^2} \sum_{i=1}^m \EE\left[\norm{a_i }^2\right], 
	\end{equation}
and for any consistent random variable $b$ being independent of $a_i$, it holds 
\begin{equation}\label{eq:dd}
\EE\left[ \norm{a_i + b}^2 \right] = \EE\left[\norm{a_i}^2 + \norm{b}^2\right].
\end{equation}
\end{lemma}

\begin{lemma}
	For any matrix $X\in\RR^{m\times d}$,  it holds that for the projection matrix $ \bPi $ defined in Eq.~\eqref{eq:Pi_def},
	\begin{equation} \label{eq:Pi}
		\norm{\bPi X} \le \norm{X}.
	\end{equation}
\end{lemma}

\begin{lemma}[Lemma 6 of \cite{qu2017harnessing}] Let Assumption~\ref{ass:f} hold, then
	\begin{align}
		\norm{\nabla f(\bbx_t) - \frac{1}{m} \mathbf{1}^\top \nabla F(\xb_t)} \le \frac{L}{\sqrt{m}} \norm{\bPi\xb_t}. \label{eq:gg}
	\end{align}
\end{lemma}

\begin{lemma}[Lemma 3 of \cite{song2021optimal}]
	Let $ f_i:\RR^d \to \RR $ satisfy Assumption~\ref{ass:f}. Denoting that $ \bbg_t = \frac{1}{m}\sum_{i=1}^m \nabla f_i(\xbti) $, it holds that
	\begin{equation}\label{eq:xx_star}
		f(\bbx_t) \le f(x^*) + \dotprod{\bbg_t, \bbx_t - x^*} - \frac{\mu}{4} \norm{\bbx_t - x^*}^2 + \frac{L}{m}\norm{\bPi\xb_t}^2.
	\end{equation}
\end{lemma}

\section{ Important Lemmas Related to Our Algorithms}

\begin{lemma}
	Letting Assumption~\ref{ass:f} hold, then we have the following inequalities:
	\begin{align}
		&\norm{\frac{1}{m} \sum_{i=1}^m \nabla f_i(\xbti) - \frac{1}{m}\sum_{i=1}^m\nabla f_i(\qbti)  }^2 
		\le \frac{1}{m} \norm{\nabla F(\xb_t) - \nabla F(\qb_t)}^2, \label{eq:gq}\\
		&\norm{\frac{1}{m}\sum_{i=1}^m \nabla f_i(\xbti)}^2 
		\le \frac{2L^2}{m}\norm{\bPi\xb_t}^2 + 4L \big( f(\bbx_t) - f(x^*) \big). \label{eq:grad_up}
	\end{align}
\end{lemma}
\begin{proof}
	For the first inequality, we have
	\begin{align*}
		\norm{\frac{1}{m} \sum_{i=1}^m\left(\nabla f_i(\xbti) - \nabla f_i(\qbti)\right) }^2
		\le \frac{1}{m}\sum_{i=1}^m \norm{\nabla f_i(\xbti) - \nabla f_i(\qbti)  }^2
		= \frac{1}{m} \norm{\nabla F(\xb_t) - \nabla F(\qb_t)}^2.
	\end{align*}
For the second inequality, we have
	\begin{align*}
		\norm{\frac{1}{m}\sum_{i=1}^m \nabla f_i(\xbti)}^2 
		=& 2\norm{\frac{1}{m}\sum_{i=1}^m \nabla f_i(\xbti) - \nabla f(\bbx_t)}^2 + 2\norm{\nabla f(\bbx_t)}^2
		\stackrel{\eqref{eq:gg}}{\le}\frac{2L^2}{m}\norm{\bPi\xb_t}^2 + 2\norm{\nabla f(\bbx_t)}^2\\
		\le& \frac{2L^2}{m}\norm{\bPi\xb_t}^2 + 4L \big( f(\bbx_t) - f(x^*) \big),
	\end{align*}
where the last inequality is because $ f(\cdot) $ is $ L $-smooth implied by Assumption~\ref{ass:f}.
\end{proof}

\begin{lemma}
	Letting Assumption~\ref{ass:f} hold, then we have the following inequalities:
	\begin{align}
		\norm{ \nabla F(\xb_t) - \nabla F(\mathbf{1} x^*) }^2 
		\le& 2L^2 \norm{\bPi \xb_t}^2 +4mL  \big( f(\bbx_t) - f(x^*) \big) \label{eq:aa}, \\
		\norm{\nabla F(\xb_t) - \nabla F(\qb_t)}^2 
		\le& 4L^2\norm{\bPi \xb_t}^2 + 2\norm{ \nabla F(\qb_t) -\nabla F(\mathbf{1} x^*) }^2 + 8mL  \big( f(\bbx_t) - f(x^*) \big). \label{eq:xq}
	\end{align}
\end{lemma}
\begin{proof}
	First, we have
	\begin{align*}
		\norm{ \nabla F(\xb_t) - \nabla F(\mathbf{1} x^*) }^2 
		\le& 2\norm{ \nabla F(\xb_t) -\nabla F(\mathbf{1}\bbx_t) }^2 + 2 \norm{ \nabla F(\mathbf{1}\bbx_t) - \nabla F(\mathbf{1} x^*)}^2\\
		\le& 2L^2 \norm{\bPi \xb_t}^2 + 2 \norm{ \nabla F(\mathbf{1}\bbx_t) - \nabla F(\mathbf{1} x^*)}^2,
	\end{align*}
where the last inequality is because of $ \nabla F(\xb) $ is $ L $-smooth implied by Assumption~\ref{ass:f}. 
In addition, applying the $ L $-smoothness of $ f_i $ again, we can obtain that
	\begin{align*}
		\norm{\nabla F(\mathbf{1}\bbx_t) - \nabla F(\mathbf{1} x^*)}^2 
		=& \sum_{i=1}^m \norm{\nabla f_i(\bbx_t) - \nabla f_i(x^*)}^2\\
		\le& \sum_{i=1}^m 2L\left( f_i(\bbx_t) - f_i(x^*) - \dotprod{\nabla f_i(x^*), \bbx_t - x^*} \right)\\
		=& 2 m L \left( f(\bbx_t) - f(x^*) - \dotprod{ \frac{1}{m} \sum_{i=1}^m\nabla f_i(x^*), \bbx_t - x^* } \right)\\
		=& 2m L\big( f(\bbx_t) - f(x^*) \big),
	\end{align*}
	where the last equality is because of $\frac{1}{m} \sum_{i=1}^m\nabla f_i(x^*) = \nabla f(x^*) = 0$. 
	Combining above two inequality, we can obtain Eq.~\eqref{eq:aa}.
	
	Furthermore, 
		\begin{align*}
		\norm{\nabla F(\xb_t) - \nabla F(\qb_t)}^2 
		\le& 2\norm{ \nabla F(\xb_t) - \nabla F(\mathbf{1} x^*) }^2  + 2\norm{ \nabla F(\qb_t) -\nabla F(\mathbf{1} x^*) }^2\\
		\stackrel{\eqref{eq:aa}}{\le}&
		4L^2\norm{\bPi \xb_t}^2 + 2\norm{ \nabla F(\qb_t) -\nabla F(\mathbf{1} x^*) }^2 + 8mL  \big( f(\bbx_t) - f(x^*) \big),
	\end{align*}
which concludes the proof.
\end{proof}

%\begin{lemma}
%	Letting Assumption~\ref{ass:noise} hold, then it holds that
%	\begin{equation}\label{eq:ex_up}
%		\EE\left[\norm{\frac{1}{m}\sum_{i=1}^m\left(\nabla f_i(\xbti,\xi_t^{(i)}) - \nabla f_i(\xbti)\right)}^2\right]
%		\le \frac{\bsig^2}{m}.
%	\end{equation}
%\end{lemma}
%\begin{proof}
%	It holds that
%	\begin{align*}
%		\EE\left[\norm{\frac{1}{m}\sum_{i=1}^m\left(\nabla f_i(\xbti,\xi_t^{(i)}) - \nabla f_i(\xbti)\right)}^2\right]
%		\stackrel{\eqref{eq:ee}}{=} \frac{1}{m^2}\sum_{i=1}^m \EE\left[\norm{\nabla f_i(\xbti,\xi_t^{(i)}) - \nabla f_i(\xbti)}^2\right]
%		\le \frac{\bsig^2}{m},
%	\end{align*}
%where the last inequality is because of Assumption~\ref{ass:noise}.
%\end{proof}

\begin{lemma}
		Letting Assumption~\ref{ass:noise} hold, then it holds that 
	\begin{equation}\label{eq:nn}
		\EE\left[ \norm{ \nabla F(\xb_t, \xi_t) - \nabla F(\qb_t,\xi_\tau)  }^2 \right]
		\le 2m\bsig^2 + \norm{ \nabla F(\xb_t) - \nabla F(\qb_t)}^2.
	\end{equation}
\end{lemma}
\begin{proof}
	First, using the fact that $\EE\left[ \nabla F(\xb_t, \xi_t) - \nabla F(\qb_t,\xi_\tau) - \left(\nabla F(\xb_t) - \nabla F(\qb_t) \right) \right] = 0$, we can obtain that
	\begin{align*}
		&\EE\left[ \norm{ \nabla F(\xb_t, \xi_t) - \nabla F(\qb_t,\xi_\tau)  }^2 \right]\\
		\stackrel{\eqref{eq:dd}}{=}&  \EE\left[\norm{\nabla F(\xb_t, \xi_t) - \nabla F(\qb_t,\xi_\tau) - \left(\nabla F(\xb_t) - \nabla F(\qb_t) \right)}^2\right] + \norm{ \nabla F(\xb_t) - \nabla F(\qb_t)}^2.
	\end{align*}
	Similarly, it holds that $\EE\left[ \nabla F(\xb_t, \xi_t) - \nabla F(\xb_t) \right] = 0$ and $\EE\left[\nabla F(\qb_t, \xi_\tau) - \nabla F(\qb_t)\right] = 0$.
	Consequently, 
	\begin{align*}
		&\EE\left[\norm{\nabla F(\xb_t, \xi_t) - \nabla F(\qb_t,\xi_\tau) - \left(\nabla F(\xb_t) - \nabla F(\qb_t) \right)}^2\right]\\
		=& \EE\left[\norm{\nabla F(\xb_t, \xi_t) - \nabla F(\xb_t)  }^2 + \norm{\nabla F(\qb_t,\xi_\tau) - \nabla F(\qb_t)}^2\right]\\
		=& \EE\left[\sum_{i=1}^m \left(\norm{\nabla f_i(\xbti,\xi_t^{(i)}) - \nabla f_i(\xbti)}^2 + \norm{\nabla f_i(\qbti, \xi_\tau^{(i)}) - \nabla f_i(\xbti)}^2\right)\right]\\
		\le& 2m\bsig^2.
	\end{align*}
	Combining above two equations, we can obtain that
	\begin{equation*}
		\EE\left[ \norm{ \nabla F(\xb_t, \xi_t) - \nabla F(\qb_t,\xi_\tau)  }^2 \right]
		\le 2m\bsig^2 + \norm{ \nabla F(\xb_t) - \nabla F(\qb_t)}^2.
	\end{equation*}
\end{proof}

\begin{lemma}\label{lem:av}
	Let $ A $, $ B $, and $ C $ be three positive constants. 
	Non-negative sequences $ \{r_t\} $, $ \{e_t\} $, and $ \{\eta_t\} $ satisfies the follow property
	\begin{align}
		r_{t+1} \le \exp\left( - \frac{\mu\eta_t}{2}\right) r_t - \eta_t e_t A + \eta_t^2 B + \eta_t^3 C. \label{eq:rr}
	\end{align}
	Set step size sequence and weight sequence as follows
	\begin{align*}
		\eta_t = 6\beta \left(L + \beta \mu t \right)^{-1},\quad \mbox{ and } \quad \omega_t = \frac{\eta_t}{\eta_0} \exp\left(\frac{\mu}{2} \sum_{i=0}^t \eta_i \right).
	\end{align*}
	Letting $S_T = \sum_{t=0}^T \omega_t$, then it holds that
	\begin{align}\label{eq:SS}
		\frac{A}{S_T} \sum_{t=0}^T e_t \omega_t 
		\le
		\frac{L^3}{2\beta^3\mu^2} \cdot \frac{1}{T^3} \cdot r_0 
		+ \frac{18 B L^3 (L + \beta\mu (T+1))^2}{\beta^2(L - \beta \mu)^3 \mu^3 T^3}  + \frac{108L^3}{\mu^2(L - \beta \mu)^3} \cdot \frac{C}{T^2} 
	\end{align}
\end{lemma}

\begin{proof}
	Diving $ \eta_t $ both sides of Eq.~\eqref{eq:rr} and rearranging it, we can obtain that
	\begin{equation}\label{eq:sum}
		\begin{aligned}
			A\sum_{t=0}^T \omega_t e_t 
			\le& 
			\sum_{t=0}^T \left(\frac{\exp\left(-\frac{\mu\eta_t}{2}\right) \omega_t }{\eta_t} \cdot r_t - \frac{\omega_t}{\eta_t} r_{t+1} + B  \omega_t\eta_t  + C\omega_t \eta_t^2\right)\\
			=& \frac{\exp\left(-\frac{\mu\eta_0}{2}\right) \omega_0 }{\eta_0} \cdot r_0 - \frac{\omega_T}{\eta_T}r_{T+1} + B \sum_{t=0}^T \omega_t\eta_t  + C\sum_{t=0}^T\omega_t \eta_t^2,
		\end{aligned}
	\end{equation}
	where the equality is because of $\omega_t = \frac{\eta_t}{\eta_{t-1}} \exp\left( \frac{\mu}{2} \eta_t \right) \omega_{t-1}$.
%	\begin{align}
%		\omega_t 
%		= \frac{\eta_t}{\eta_0} \exp\left(\frac{\mu}{2} \sum_{i=0}^t \eta_i \right) \omega_0 
%		= \frac{L}{L + \beta\mu t}\exp\left(\frac{\mu}{2} \sum_{i=0}^t \eta_i \right).
%	\end{align}

Now, we are going to upper bound $\omega_t$ as follows
	%\begin{align}
	%\omega_t 
	%=& \frac{L}{L + \beta\mu t}\exp\left(3\beta \mu \sum_{i=0}^t (L + \beta \mu i)^{-1} \right)\\
	%\le& \frac{L}{L + \beta\mu t}\exp\left(3\beta \mu \int_{0}^{t+1}  (L + \beta \mu i)^{-1}\; di\right)\\
	%\le& \frac{L}{L + \beta\mu t} \cdot \left( \frac{L + \beta\mu(t+1)}{L}-1\right)^3\\
	%=& \frac{L}{L + \beta\mu t} \left(\frac{\beta\mu(t+1)}{L}\right)^3
	%\end{align}
	\begin{equation}\label{eq:w_up}
		\begin{aligned}
			\omega_t 
			=& \frac{L}{L + \beta\mu t}\exp\left(3\beta \mu \sum_{i=0}^t (L + \beta \mu i)^{-1} \right)
			\stackrel{\eqref{eq:int_2}}{\le} \frac{L}{L + \beta\mu t}\exp\left(3\beta \mu \int_{-1}^{t}  (L + \beta \mu i)^{-1}\; di\right)\\
			=& \frac{L}{L + \beta\mu t} \cdot \left( \frac{L + \beta\mu t}{L - \beta \mu}\right)^3
			= \frac{L(L + \beta\mu t)^2}{(L - \beta \mu)^3}.
		\end{aligned}
	\end{equation}
	We lower bound $ S_T$ as follows:
	\begin{equation}\label{eq:S_T}
		\begin{aligned}
			S_T 
			=& \sum_{t=0}^T \omega_t 
			= \sum_{t=0}^T \frac{L}{L + \beta\mu t}\exp\left(3\mu\beta \sum_{i=0}^t (L+\beta\mu i)^{-1}\right)\\
			\stackrel{\eqref{eq:int_2}}{\ge}& 
			\sum_{t=0}^T \frac{L}{L + \beta\mu t}\exp\left(3\mu\beta \int_{i=0}^{t} (L+\beta\mu i)^{-1} \;di\right)\\
			=& \sum_{t=0}^T \frac{L}{L + \beta\mu t} \left( \frac{L + \beta\mu t}{L}\right)^3
			\stackrel{\eqref{eq:int_1}}{\ge} \int_{0}^T \left( \frac{L + \beta\mu t}{L}\right)^2 \;dt\\
			=& \frac{L}{3\beta \mu} \left( \left( 1 + \frac{\beta \mu}{L} T \right)^3 - 1 \right)
			\ge \frac{\beta^2 \mu^2}{3L^2} T^3.
		\end{aligned}	
	\end{equation}
	We also have
	\begin{equation}\label{eq:w_eta}
		\begin{aligned}
			\sum_{t=0}^T \omega_t \eta_t
			\stackrel{\eqref{eq:w_up}}{\le}&
			\sum_{t=0}^T \frac{L(L + \beta\mu t)^2}{(L - \beta \mu)^3} \cdot \frac{6\beta}{L + \beta \mu t}
			\stackrel{\eqref{eq:int_1}}{\le}
			\frac{6L\beta}{(L - \beta \mu)^3} \int_{0}^{T+1} (L+\beta \mu t)\; dt \\
			\le&
			\frac{6L (L + \beta\mu (T+1))^2}{(L - \beta \mu)^3 \mu},
		\end{aligned}	
	\end{equation}
and
	\begin{equation}\label{eq:w_eta_2}
		\begin{aligned}
			\sum_{t=0}^T\omega_t \eta_t^2
			\stackrel{\eqref{eq:w_up}}{\le}
			\sum_{t=0}^T \frac{L(L + \beta\mu t)^2}{(L - \beta \mu)^3} \cdot \left(\frac{6\beta}{L + \beta \mu t}\right)^2
			=
			\frac{36L\beta^2}{(L - \beta \mu)^3} T.
		\end{aligned}	
	\end{equation}

	Dividing $S_T$ both sides of Eq.~\eqref{eq:sum}, we can obtain that
	\begin{align*}
		\frac{A}{S_T} \sum_{i=0}^T \omega_t e_t
		\le&
		\frac{\exp\left(-\frac{\mu\eta_0}{2}\right) \omega_0 }{S_T\eta_0} \cdot r_0  + \frac{B}{S_T} \sum_{t=0}^T \omega_t\eta_t  + \frac{C}{S_T}\sum_{t=0}^T\omega_t \eta_t^2 \\
		\stackrel{\eqref{eq:S_T}\eqref{eq:w_eta}\eqref{eq:w_eta_2}}{\le}&
		\frac{L^3}{2\beta^3\mu^2} \cdot \frac{1}{T^3} \cdot r_0 
		+ \frac{18 B L^3 (L + \beta\mu (T+1))^2}{\beta^2(L - \beta \mu)^3 \mu^3 T^3}  + \frac{108L^3}{\mu^2(L - \beta \mu)^3} \cdot \frac{C}{T^2}.
	\end{align*}
\end{proof}

\section{Proofs for Section~\ref{sec:ss}}

\subsection{Proof of Lemma~\ref{lem:s}}
\begin{proof}[Proof of Lemma~\ref{lem:s}]
	We prove the result by the induction.
	For $t = 0$, by the initialization, we have
	\begin{equation*}
		\bbs_0 = \frac{1}{m} \sum_{i=1}^m \nabla f_i(\qb_t^{(0)}, \xi_0^{(i)}).
	\end{equation*}
	If $\zeta_t = 0$, then $\qb_{t+1} = \qb_t$. By the induction hypothesis
	and the fact $\mathbf{1}^\top W = \mathbf{1}^\top$, we have
	\begin{equation*}
		\bbs_{t+1} = \bbs_t = \frac{1}{m}\sum_{i=1}^m \nabla f_i(\qbti, \xi_\tau^{(i)}) = \frac{1}{m}\sum_{i=1}^m \nabla f_i(\qb_{t+1}^{(i)}, \xi_\tau^{(i)}).
	\end{equation*}
	If $\zeta_t = 1$, then $\qb_{t+1} = \xb_t$. By the induction hypothesis, 
	\begin{align*}
		\bbs_{t+1} 
		=& \bbs_t + \frac{1}{m}\sum_{i=1}^m \left(\nabla f_i(\xbti,\xi_t^{(i)}) - \nabla f_i(\qbti,\xi_\tau^{(i)})\right)\\
		=& \frac{1}{m}\sum_{i=1}^m \nabla f_i(\qbti,\xi_\tau^{(i)}) +  \frac{1}{m}\sum_{i=1}^m \left(\nabla f_i(\xbti,\xi_t^{(i)}) - \nabla f_i(\qbti,\xi_\tau^{(i)})\right) \\
		=& \frac{1}{m}\sum_{i=1}^m\nabla f_i(\xbti,\xi_t^{(i)}) 
		=  \frac{1}{m}\sum_{i=1}^m\nabla f_i(\qb_{t+1}^{(i)},\xi_t^{(i)}) 
		= \frac{1}{m}\sum_{i=1}^m\nabla f_i(\qb_{t+1}^{(i)},\xi_\tau^{(i)}),
	\end{align*}
	where the last equality is because of $\tau = t$ if $\zeta_t = 1$.
\end{proof}

\subsection{Proof of Lemma~\ref{lem:main_1}}

Before the detailed proof, we first introduce the following lemma which describes the evolution of consensus error terms.
\begin{lemma}\label{lem:err_consensus}
	Suppose Assumptions~\ref{ass:noise}-\ref{ass:mix} hold. Sequences $ \{\xb_t\} $, $ \{\bs_t\} $, and $ \{\qb_t\} $ are generated by Algorithm~\ref{alg:alg_name}. We have the following inequalities: 
	\begin{equation}\label{eq:xx}
		\begin{aligned}
			&\EE\left[\norm{\bPi \xb_{t+1}}^2 \mid \cF_t\right]\\
			\le& \left(1 - \theta\right) \cdot \norm{\bPi\xb_t}^2 + \frac{3\eta_t^2}{\theta} \norm{\bPi \bs_t}^2 +  \left( {\eta_t^2} + \frac{2\eta_t^2}{\theta} \right) \norm{ \nabla F(\xb_t) - \nabla F(\qb_t)}^2 + {2m\eta_t^2\bsig^2}
		\end{aligned}
	\end{equation}
	and
	\begin{equation}\label{eq:ss}
		\EE\left[\norm{\bPi\bs_{t+1}}^2 \mid \cF_t\right]
		\le \left(1 - \theta\right) \norm{\bPi\bs_t}^2 
		+ 2p  \norm{\nabla F(\xb_t) - \nabla F(\qb_t)}^2 + 4mp\bsig^2
	\end{equation}
	and
	\begin{equation}\label{eq:qx}
		\begin{aligned}
			& \EE\left[\norm{\nabla F(\qb_{t+1}) - \nabla F(\mathbf{1}x^*)}^2 \mid \cF_t \right]\\
			\le& (1 -p) \norm{ \nabla F(\qb_t) - \nabla F(\mathbf{1}x^*) }^2 + 2pL^2\norm{\bPi \xb_t}^2  + 4mLp \big(f(\bbx_t) - f(x^*)\big). 
		\end{aligned}
	\end{equation}
\end{lemma}
\begin{proof}
	First, we have
	\begin{equation}\label{eq:x_1}
		\begin{aligned}
			&\EE\left[\norm{\bPi\left( \xb_t - \eta_t \left( \bs_t +  \nabla F(\xb_t, \xi_t) - \nabla F(\qb_t,\xi_\tau)  \right) \right)  }^2\mid \cF_t\right]	\\
			=& \EE\left[ \norm{ \bPi\left( \xb_t - \eta_t \bs_t \right) - \eta_t \cdot \bPi\left( \nabla F(\xb_t, \xi_t) - \nabla F(\qb_t,\xi_\tau) \right) }^2 \mid \cF_t \right]\\
			=&  \norm{\bPi(\xb_t -\eta_t \bs_t)}^2 + \eta_t^2 \cdot \EE\left[\norm{ \bPi \big( \nabla F(\xb_t, \xi_t) - \nabla F(\qb_t,\xi_\tau) \big)  }^2 \mid \cF_t\right]\\
			& -2 \eta_t \cdot \EE \left[\dotprod{\bPi(\xb_t - \eta_t \bs_t), \bPi \big( \nabla F(\xb_t, \xi_t) - \nabla F(\qb_t,\xi_\tau) \big) } \mid \cF_t\right]\\
			=& \norm{\bPi(\xb_t -\eta_t \bs_t)}^2 + \eta_t^2 \cdot  \EE\left[ \norm{ \bPi \big( \nabla F(\xb_t, \xi_t) - \nabla F(\qb_t,\xi_\tau) \big)  }^2 \mid \cF_t\right] \\
			& -2\eta_t\dotprod{\bPi(\xb_t - \eta_t \bs_t),   \bPi \big( \nabla F(\xb_t) - \nabla F(\qb_t) \big)}\\
			\le& \norm{\bPi(\xb_t -\eta_t \bs_t)}^2 + \eta_t^2  \EE\left[ \norm{ \bPi \big( \nabla F(\xb_t, \xi_t) - \nabla F(\qb_t,\xi_\tau) \big)  }^2 \mid \cF_t \right]\\
			& + \frac{\theta}{2} \norm{ \bPi (\xb_t -\eta_t \bs_t) }^2 + \frac{2\eta_t^2}{\theta} \norm{ \bPi \big( \nabla F(\xb_t) - \nabla F(\qb_t) \big)  }^2\\
			\stackrel{\eqref{eq:Pi}}{\le}& 
			\left(1 + \frac{\theta}{2}\right) \norm{\bPi(\xb_t -\eta_t \bs_t)}^2 + \eta_t^2  \EE\left[ \norm{ \nabla F(\xb_t, \xi_t) - \nabla F(\qb_t,\xi_\tau)  }^2 \right] \\
			& + \frac{2\eta_t^2}{\theta} \norm{  \nabla F(\xb_t) - \nabla F(\qb_t) }^2\\
			\stackrel{\eqref{eq:nn}}{\le}&
			\left(1 + \frac{\theta}{2}\right) \norm{\bPi(\xb_t -\eta_t \bs_t)}^2 
			+ \left( {\eta_t^2} + \frac{2\eta_t^2}{\theta} \right)\norm{  \nabla F(\xb_t) - \nabla F(\qb_t) }^2 + {2m\eta_t^2\bsig^2},
		\end{aligned}	
	\end{equation}
	where the first inequality is because of the Cauchy's inequality. 
	
	Using above equation, we can obtain that
	\begin{align*}
		&\EE\left[\norm{\bPi \xb_{t+1}}^2 \mid \cF_t\right]\\
		=& \EE\left[\norm{\bPi W \left( \xb_t - \eta_t \left( \bs_t + \zeta_t \left( \nabla F(\xb_t,\xi_t) - \nabla F(\qb_t, \xi_\tau) \right) \right) \right)  }^2 \mid \cF_t\right]\\
		\stackrel{\eqref{eq:dec}}{\le}&
		\left(1-\theta\right)^2 \cdot  \EE\left[\norm{\bPi\left( \xb_t - \eta_t \left( \bs_t + \zeta_t \left( \nabla F(\xb_t, \xi_t) - \nabla F(\qb_t,\xi_\tau) \right) \right) \right)  }^2 \mid \cF_t\right]\\
		\stackrel{\eqref{eq:x_1}}{\le}&
		\left(1-\theta\right)^2 \left(1 + \frac{\theta}{2}\right)\norm{\bPi(\xb_t -\eta_t \bs_t)}^2 + \left( {\eta_t^2} + \frac{2\eta_t^2}{\theta} \right)\norm{  \nabla F(\xb_t) - \nabla F(\qb_t) }^2 + {2m\eta_t^2\bsig^2} \\
		\le& \left(1-\theta\right)^2 \left(1 + \frac{\theta}{2}\right) \left( \left( 1+\frac{\theta}{2} \right) \norm{\bPi\xb_t} + \left(1 + \frac{2}{\theta}\right)\eta_t^2 \norm{\bPi\bs_t}^2 \right)\\
		&+\left( {\eta_t^2} + \frac{2\eta_t^2}{\theta} \right)\norm{  \nabla F(\xb_t) - \nabla F(\qb_t) }^2 + {2m\eta_t^2\bsig^2}\\
		\le& 
		\left(1 - \theta\right) \norm{\bPi\xb_t}^2 + \frac{3\eta_t^2}{\theta}\norm{\bPi\bs_t}^2+ \frac{3\eta_t^2}{\theta} \norm{  \nabla F(\xb_t) - \nabla F(\qb_t) }^2 + {2m\eta_t^2\bsig^2},
	\end{align*}
	which proves the result of Eq.~\eqref{eq:xx}.
	
	Now, we are going to prove Eq.~\eqref{eq:ss} and we have
	\begin{align*}
		&\EE\left[\norm{\bPi\bs_{t+1}}^2 \mid \cF_t\right]\\
		=& (1 - p) \EE\left[\norm{\bPi W \bs_t}^2 \mid \cF_t \right]+ p  \EE\left[\norm{\bPi(W\bs_t + \nabla F(\xb_t,\xi_t) - \nabla F(\qb_t,\xi_\tau))}^2 \mid \cF_t\right]\\
		\le& (1 - p) \EE\left[\norm{\bPi W \bs_t}^2 \mid \cF_t\right] + 2p \EE\left[  \norm{\bPi W\bs_t}^2 +  \norm{\nabla F(\xb_t,\xi_t) - \nabla F(\qb_t,\xi_\tau)}^2 \mid \cF_t \right]\\
		=& (1 + p) \EE\left[\norm{\bPi W\bs_t}^2 \mid \cF_t\right] + 2p \EE\left[ \norm{\nabla F(\xb_t,\xi_t) - \nabla F(\qb_t,\xi_\tau)}^2 \mid \cF_t\right] \\ 
		\stackrel{\eqref{eq:dec}}{\le}& 
		(1+p)\left(1 - \theta\right)^2 \cdot \norm{\bPi\bs_t}^2 + 2p \EE\left[ \norm{\nabla F(\xb_t,\xi_t) - \nabla F(\qb_t,\xi_\tau)}^2 \mid \cF_t\right]\\
		\le& \left(1 - \theta\right) \EE\left[\norm{\bPi\bs_t}^2 \right]
		+ 2p \EE\left[\norm{\nabla F(\xb_t,\xi_t) - \nabla F(\qb_t,\xi_\tau)}^2 \mid \cF_t\right],
	\end{align*}
	where the last inequality is because of $p \le \theta$.
	Using Eq.~\eqref{eq:nn}, we can obtain that
	\begin{align*}
		\EE\left[\norm{\bPi\bs_{t+1}}^2 \mid \cF_t\right]
		\stackrel{\eqref{eq:nn}}{\le}
		\left(1 - \theta\right)\norm{\bPi\bs_t}^2 
		+ 2p \norm{\nabla F(\xb_t) - \nabla F(\qb_t)}^2 + 4mp\bsig^2.
	\end{align*}
	
	By the update rule of $ \qb_t $, we have
	\begin{align*}
		&\EE\left[\norm{\nabla F(\qb_{t+1}) - \nabla F(\mathbf{1}x^*)}^2 \mid \cF_t\right] \\
		\stackrel{\eqref{eq:qq}}{=}& p \norm{\nabla F(\xb_t) - \nabla F(\mathbf{1}x^*)}^2 + (1 -p) \norm{ \nabla F(\qb_t) - \nabla F(\mathbf{1}x^*) }^2 \\
		\stackrel{\eqref{eq:aa}}{\le}& 2pL^2\norm{\bPi \xb_t}^2  + 4mLp \big(f(\bbx_t) - f(x^*)\big) + (1 -p) \norm{ \nabla F(\qb_t) - \nabla F(\mathbf{1}x^*) }^2.
	\end{align*}
\end{proof}

By the above lemma and the setting of parameters, we can prove Lemma~\ref{lem:main_1} as follows.
\begin{proof}[Proof of Lemma~\ref{lem:main_1}] First, we have
	\begin{align*}
		&\EE\left[\norm{\bPi \xb_{t+1}}^2 + C_{1,t+1} \norm{\bPi\bs_{t+1}}^2 + C_{2,t+1} \norm{\nabla F(\qb_{t+1}) - \nabla F(\mathbf{1}x^*)}^2\right]\\
		\stackrel{\eqref{eq:xx}\eqref{eq:ss}\eqref{eq:qx}}{\le}& 
		\left(1 - \theta + 2pL^2 C_{2,t}
		\right)\cdot \norm{\bPi\xb_t}^2
		+\left(1 - \theta + \frac{3\eta_t^2}{\theta C_{1,t}}\right) C_{1,t} \norm{\bPi\bs_t}^2\\
		&+\left(1 - p\right)C_{2, t}\norm{\nabla F(\qb_t) - \nabla F(\mathbf{1}x^*)}^2 
		+ \left({\eta_t^2} +\frac{2\eta_t^2}{\theta} + {2pC_{1,t}}\right) \norm{\nabla F(\xb_t) - \nabla F(\qb_t)}^2\\
		& +\left({2m\eta_t^2} + {4mp C_{1,t}}\right) \bsig^2
		+ 4mLp C_{2,t}\big(f(\bbx_t) - f(x^*)\big)\\
		\stackrel{\eqref{eq:xq}}{\le}&
		\left(1 - \theta + 2pL^2 C_{2,t} + {4L^2\eta_t^2} + \frac{8L^2\eta_t^2}{\theta} + 8pL^2 C_{1,t}
		\right)\cdot \norm{\bPi\xb_t}^2 \\
		&+\left(1 - \theta + \frac{3\eta_t^2}{\theta C_{1,t}}\right) C_{1,t} \norm{\bPi\bs_t}^2\\
		&+\left(1 - p+ \frac{2\eta_t^2}{C_{2,t}} +\frac{4\eta_t^2}{\theta C_{2,t}} + 
		\frac{4pC_{1,t}}{C_{2,t}}\right)C_{2, t}\norm{\nabla F(\qb_t) - \nabla F(\mathbf{1}x^*)}^2 \\
		& +\left({2m\eta_t^2} + {4mp C_{1,t}}\right) \bsig^2
		+ 4mL \left(p C_{2,t} + {2\eta_t^2} +\frac{4\eta_t^2}{\theta} + {4pC_{1,t}}\right)\big(f(\bbx_t) - f(x^*)\big)\\
		=& 
		\left(1 - \theta + 4L\eta_t + \frac{( 4\theta + 8 + 32) L^2 \eta_t^2 }{\theta}\right) \cdot \norm{\bPi\xb_t}^2 + \left(1 -\theta + \frac{3\theta}{4}\right) C_{1,t} \norm{\bPi\bs_t}^2\\
		& + \left(1 - \theta + \eta_t L\theta + 2L\eta_t +  8\eta_t L\right) \cdot C_{2, t}\norm{\nabla F(\qb_t) - \nabla F(\mathbf{1}x^*)}^2\\
		& + \left( {2m\eta_t^2} + \frac{16m\eta_t^2}{\theta} \right)\bsig^2 +  4mL \left( \frac{2\eta_t}{L} + {2\eta_t^2} + \frac{4\eta_t^2}{\theta} + \frac{16\eta_t^2}{\theta} \right)\big(f(\bbx_t) - f(x^*)\big)\\
		\le&
		\left(1 - \frac{\theta}{4}\right)\cdot \left( \norm{\bPi\xb_t}^2 + C_{1,t} \norm{\bPi\bs_t}^2 +  C_{2, t}\norm{\nabla F(\qb_t) - \nabla F(\mathbf{1}x^*)}^2 \right)\\
		& +  \frac{18m\eta_t^2}{\theta} \bsig^2 + \left( \frac{88mL\eta_t^2}{\theta} +8m\eta_t \right)\big(f(\bbx_t) - f(x^*)\big),
	\end{align*}
	where the first equality is because of $C_{1,t} = 4\eta_t^2 / \theta^2$, $C_{2,t} = 2\eta_t/(L\theta)$, $p = \theta$ and the last inequality is because of $\eta_t \le \frac{\theta}{16L}$.
\end{proof}

\subsection{Proof of Lemma~\ref{lem:main_2}}

\begin{proof}[Proof of Lemma~\ref{lem:main_2}] By the update rule of $ \xb_t $, we have
	\begin{equation}\label{eq:xx_1}
		\begin{aligned}
			&\EE\left[ \norm{\bbx_{t+1} - x^*}^2 \mid \cF_t \right]\\
			=&\norm{ \bbx_t - x^* }^2 - 2\eta_t \EE\left[\dotprod{\bbs_t +  \frac{1}{m}\sum_{i=1}^m\left( \nabla f_i(\xbti, \xi_t^{(i)}) - \nabla f_i(\qbti),\xi_\tau^{(i)}\right), \bbx_t - x^*} \mid \cF_t\right]\\
			&+\eta_t^2\EE\left[ \norm{\bbs_t +  \frac{1}{m}\sum_{i=1}^m\left( \nabla f_i(\xbti, \xi_t^{(i)}) - \nabla f_i(\qbti),\xi_\tau^{(i)}\right)^2  }  \mid \cF_t\right]\\
			\stackrel{\eqref{eq:s}}{=}&
			\norm{ \bbx_t - x^* }^2 -2 \eta_t \EE\left[\dotprod{\frac{1}{m}\sum_{i=1}^m \nabla f_i(\xbti, \xi_t^{(i)}), \bbx_t - x^*} \mid \cF_t\right] \\
			&+ \eta_t^2\EE\left[ \norm{\frac{1}{m}\sum_{i=1}^m\nabla f_i(\xbti, \xi_t^{(i)}) }^2 \mid \cF_t  \right]\\
			=&
			\norm{ \bbx_t - x^* }^2 -2 \eta_t \dotprod{\frac{1}{m}\sum_{i=1}^m \nabla f_i(\xbti), \bbx_t - x^*} + \eta_t^2\EE\left[ \norm{\frac{1}{m}\sum_{i=1}^m\nabla f_i(\xbti, \xi_t^{(i)}) }^2 \mid \cF_t  \right] \\
			\stackrel{\eqref{eq:xx_star}}{\le}& 
			\left(1 - \frac{\mu\eta_t}{2}\right) \norm{\bbx_t - x^*}^2 - 2\eta_t\big( f(\bbx_t) - f(x^*) \big) + \frac{2L\eta_t}{m} \norm{\bPi\xb_t}^2 \\
			&+ \eta_t^2\EE\left[ \norm{\frac{1}{m}\sum_{i=1}^m\nabla f_i(\xbti, \xi_t^{(i)}) }^2 \mid \cF_t  \right].
		\end{aligned}
	\end{equation}
	
	Furthermore, using the fact that $ \EE\left[\nabla f_i(\xbti,\xi_t^{(i)}) - \nabla f_i(\xbti) \mid \cF_t\right] = 0$, we have
	\begin{align*}
		&\EE\left[\norm{\frac{1}{m}\sum_{i=1}^m\nabla f_i(\xbti,\xi_t^{(i)}) }^2 \mid \cF_t\right]\\
		\stackrel{\eqref{eq:dd}}{\le}& 
		\EE\left[\norm{\frac{1}{m}\sum_{i=1}^m\nabla f_i(\xbti,\xi_t^{(i)}) -\frac{1}{m} \sum_{i=1}^m\nabla f_i(\xbti)}^2 \mid \cF_t\right] + \norm{\frac{1}{m}\sum_{i=1}^m \nabla f_i(\xbti)}^2\\
		\stackrel{\eqref{eq:ee}}{=}& 
		\frac{1}{m^2} \sum_{i=1}^m\EE\left[\norm{\nabla f_i(\xbti, \xi_t^{(i)}) - \nabla f_i(\xbti)}^2 \mid \cF_t\right] 
		+ \norm{\frac{1}{m}\sum_{i=1}^m \nabla f_i(\xbti)}^2\\
		\stackrel{\eqref{eq:noise}}{\le}&
		\frac{\bsig^2}{m} + \norm{\frac{1}{m}\sum_{i=1}^m \nabla f_i(\xbti)}^2\\
		\stackrel{\eqref{eq:grad_up}}{\le}& 
		\frac{2L^2}{m}\norm{\bPi\xb_t}^2 + 4L\big(f(\bbx_t) - f(x^*)\big) + \frac{\bsig^2}{m}.
	\end{align*}
	
	Combining above results, we can obtain that 
	\begin{align*}
		&\EE\left[ \norm{\bbx_{t+1} - x^*}^2 \mid \cF_t \right]\\
		\le& \left(1 - \frac{\mu\eta_t}{2}\right) \norm{\bbx_t - x^*}^2 - 2\eta_t\big( f(\bbx_t) - f(x^*) \big) + \frac{2L\eta_t}{m} \norm{\bPi\xb_t}^2 
		+ \eta_t^2  \cdot \frac{\bsig^2}{m}\\
		&+ \frac{2\eta_t^2L^2}{m} \norm{\bPi \xb_t}^2 
		+ 4\eta_t^2 L \big(f(\bbx_t) - f(x^*)\big)\\
		=&\left(1 - \frac{\mu\eta_t}{2}\right) \norm{\bbx_t - x^*}^2 - 2\eta_t \left(1 - 2\eta_t L \right)\big(f(\bbx_t) - f(x^*)\big) + \eta_t^2\cdot \frac{\bsig^2}{m}\\
		&+\frac{2L\eta_t \left( 1 + 2\eta_t L \right) }{m} \norm{\bPi \xb_t}^2.
	\end{align*}

\end{proof}

\subsection{Proof of Lemma~\ref{lem:main_3}}

\begin{proof}[Proof of Lemma~\ref{lem:main_3}] We have
	\begin{align*}
		&\EE\left[\norm{\bbx_{t+1} - x^*}^2 + \frac{24L\eta_{t+1}}{m\theta} \Psi_{t+1} \mid \cF_t\right]\\
		\stackrel{\eqref{eq:xx_a}\eqref{eq:psi}}{\le}&
		\left(1 - \frac{\mu\eta_t}{2}\right) \norm{\bbx_t - x^*}^2 - 2\eta_t \left(1 - 2\eta_t L \right)\big(f(\bbx_t) - f(x^*)\big) + \eta_t^2\cdot \frac{\bsig^2}{m}\\
		&+\frac{2L\eta_t \left( 1 + 2\eta_t L \right) }{m} \norm{\bPi \xb_t}^2 
		+ \left(1 - \frac{\theta}{4}\right) \cdot \frac{24L\eta_t}{m\theta} \Psi_t \\
		&+ \frac{2^4 \cdot 3^3 \cdot L\eta_t^3}{\theta^2}\bsig^2 + \left(\frac{2^8\cdot 3^2\cdot L^2\eta_t^3}{\theta^2} + \frac{2^6\cdot 3\cdot L\eta_t^2}{\theta}\right)\cdot \big(f(\bbx_t) - f(x^*)\big)\\
		=&  
		\left(1 - \frac{\mu\eta_t}{2}\right) \norm{\bbx_t - x^*}^2 - 2\eta_t \left(1 - 2\eta_t L-\frac{2^7\cdot 3^2 \cdot L^2 \eta_t^2}{\theta^2} - \frac{2^5\cdot 3\cdot L\eta_t}{\theta}\right)\big(f(\bbx_t) - f(x^*)\big) \\
		&+ \eta_t^2\cdot \frac{\bsig^2}{m} + \frac{2^4 \cdot 3^3 \cdot L\eta_t^3}{\theta^2}\bsig^2 
		+\frac{3L\eta_t}{m} \norm{\bPi \xb_t}^2 
		+ \left(1 - \frac{\theta}{4}\right) \cdot \frac{24L\eta_t}{m\theta} \Psi_t \\
		\stackrel{\eta_t \le \frac{\theta}{2^6\cdot 3\cdot L}}{\le}&  
		\left(1 - \frac{\mu\eta_t}{2}\right) \norm{\bbx_t - x^*}^2 - \frac{7\eta_t}{8}\big(f(\bbx_t) - f(x^*)\big) 
		+\left( 1 - \frac{\theta}{4} + \frac{3L\eta_t}{m} \cdot \frac{m\theta}{24L\eta_t} \right) \cdot \frac{24L\eta_t}{m\theta} \Psi_t\\
		&+ \frac{2^4 \cdot 3^3 \cdot L\eta_t^3}{\theta^2}\bsig^2 + \eta_t^2\cdot \frac{\bsig^2}{m}\\
		=& 
		\left(1 - \frac{\mu\eta_t}{2}\right) \norm{\bbx_t - x^*}^2 + \left(1 - \frac{\theta}{8}\right) \cdot \frac{24L\eta_t}{m\theta} \Psi_t - \frac{7\eta_t}{8} \big(f(\bbx_t) - f(x^*)\big)\\
		& +  \frac{2^4 \cdot 3^3 \cdot L\eta_t^3}{\theta^2}\bsig^2 + \eta_t^2\cdot \frac{\bsig^2}{m}\\
		\le& 
		\left(1 - \frac{\mu\eta_t}{2}\right) \left(\norm{\bbx_t - x^*}^2 +  \frac{24L\eta_t}{m\theta} \Psi_t \right) 
		- \frac{7\eta_t}{8} \big(f(\bbx_t) - f(x^*)\big) 
		+ \frac{2^4 \cdot 3^3 \cdot L\eta_t^3}{\theta^2}\bsig^2 + \eta_t^2\cdot \frac{\bsig^2}{m}\\
		\le& \exp\left(- \frac{\mu\eta_t}{2}\right) \left(\norm{\bbx_t - x^*}^2 +  \frac{24L\eta_t}{m\theta} \Psi_t \right) 
		- \frac{7\eta_t}{8} \big(f(\bbx_t) - f(x^*)\big) 
		+ \frac{2^4 \cdot 3^3 \cdot L\eta_t^3}{\theta^2}\bsig^2 + \eta_t^2\cdot \frac{\bsig^2}{m},
	\end{align*}
	where the last inequality is because of $1 - x \le \exp(-x)$ when $0<x<1$.
\end{proof}

\section{Proofs of Section~\ref{sec:ass}}
\subsection{Proof of Lemma~\ref{lem:equal}}
\begin{proof}[Proof of Lemma~\ref{lem:equal}]
	We prove the result by the induction.
For the case $t = 0$, the it holds that $ \sum_{i=1}^m \txb_0^{(i)} =   \sum_{i=m+1}^{2m} \txb_0^{(i)}$ trivially by the initialization of Algorithm~\ref{alg:alg_namea}. 
Assuming that Eq.~\eqref{eq:equal} holds for $t$, then we have
\begin{align*}
	\sum_{i=1}^{m} \txb_{t+1}^{(i)} 
	=&
	 (1+\gamma) \mathbf{1}^\top W \txb_t^{(1:m)} - \gamma \mathbf{1}^\top W\txb_t^{(m+1 : 2m)} + \mathbf{1}^\top \left( \nabla F(\xb_t, \xi_t) - \nabla F(\qb_t, \xi_\tau)\right) \\
	=& \mathbf{1}^\top W \txb_t^{(1:m)} +  \mathbf{1}^\top \left( \nabla F(\xb_t, \xi_t) - \nabla F(\qb_t, \xi_\tau)\right)\\
	=& \sum_{i=m+1}^{2m} \txb_{t+1}^{(i)} 
\end{align*}
where the second equality is because of the induction assumption, the first and third equality are because of definition of $ \tW $ and update rule of $ \txb_t $.
The result  $\sum_{i=1}^{m} \tsb_t^{(i)} = \sum_{i=m+1}^{2m} \tsb_t^{(i)} $ can be  proved similarly.

Since $\sum_{i=1}^{m} \tsb_t^{(i)} = \sum_{i=m+1}^{2m} \tsb_t^{(i)} $, then Eq.~\eqref{eq:sq} can be proved as the same to the one of Eq.~\eqref{eq:s}.

\end{proof}

\subsection{Proof of Lemma~\ref{lem:tsx}}
\begin{proof}[Proof of Lemma~\ref{lem:tsx}]
	First, for notation convenience, we denote that
	\begin{align*}
		A_\# := \begin{bmatrix}
			A \\
			A
		\end{bmatrix} , \quad \forall \; A \in\RR^{m \times d}.
	\end{align*}
	By the update rule of $ \tsb_t $, we can obtain that
	\begin{align*}
		&\EE\left[ \norm{\tPi\tsb_{t+1}}^2 \right]\\
		=&    
		(1 - p) \norm{\tPi\tW\tPi\tsb_t}^2 + p \EE\left[\norm{ \tPi \tW \tPi \tsb_t + \tPi\bigl( \nabla F(\xb_t,\xi_t)_{\#} -\nabla F(\qb_t, \xi_\tau)_{\#}\bigr) }^2\right]\\
		\le&  
		(1+p) \norm{\tPi\tW\tPi\tsb_t}^2 + 2p \EE\left[\norm{ \nabla F(\xb_t, \xi_t)_{\#} -\nabla F(\qb_t, \xi_\tau)_{\#} }^2\right].
	\end{align*}
	Similar to above equation, we have
	\begin{align*}
		\EE  \left[\norm{\tPi\tW\tPi\tsb_t}^2\right] 
		\le (1+p) \norm{\tPi\tW^2\tPi\tsb_{t-1}  }^2 + 2p \norm{ \tPi \tW \tPi\bigl( \nabla F(\xb_{t-1}, \xi_{t-1})_{\#} -\nabla F(\qb_{t-1}, \xi_\tau)_{\#} \bigr) }^2.
	\end{align*}
	Using above equation recursively, we can obtain that
	\begin{align*}
		&\EE\left[ \norm{\tPi\tsb_{t+1}}^2 \right]
		\\
		\le& 
		(1+p) \norm{\tPi\tW\tPi\tsb_t}^2 + 2p \EE\left[\norm{ \nabla F(\xb_t, \xi_t)_{\#} -\nabla F(\qb_t, \xi_\tau)_{\#} }^2\right]\\
		\le& 
		\EE\left[ \sum_{i=0}^{t} 2p (1+p)^{t-i} \norm{ \tPi\tW^{t-i} \tPi  \bigl( \nabla F(\xb_i, \xi_i)_{\#} -\nabla F(\qb_i, \xi_{\tau_i})_{\#} \bigr) }^2  \right]\\
		&+ (1+p)^{t+1} \norm{ \tPi \tW^{t+1} \tPi \tsb_{\#,0} }^2\\
		\stackrel{\eqref{eq:tW_dec}}{\le}&
		\EE\left[ \sum_{i=0}^{t} 2\alpha p(1+p)^{t-i} \left( 1 - \tilde{\theta} \right)^{2(t-i)} \norm{ \nabla F(\xb_i, \xi_i)_{\#} -\nabla F(\qb_i, \xi_{\tau_i})_{\#} }^2 \right]\\
		& + \alpha (1+p)^{t+1} (1 - \tilde{\theta})^{2(t+1)} \norm{\tPi \tsb_{0}  }^2\\
		\le&
		\EE\left[ \sum_{i=0}^{t} 4\alpha p \left( 1 - \tilde{\theta} \right)^{t-i} \norm{ \nabla F(\xb_i,\xi_i) -\nabla F(\qb_i, \xi_{\tau_i}) }^2  \right] + 2\alpha \left(1-\tilde{\theta}\right)^{t+1} \norm{\bPi \bs_0}^2\\
		\stackrel{\eqref{eq:nn}}{\le}&
		\EE\left[ \sum_{i=0}^{t} 4\alpha p \left( 1 - \tilde{\theta} \right)^{t-i} \left( \norm{ \nabla F(\xb_i) -\nabla F(\qb_i) }^2 + 2m\bsig^2 \right) \right] + 2\alpha \left(1-\tilde{\theta}\right)^{t+1} \norm{\bPi \bs_0}^2  \\
		=& \cC_{s,t+1}
	\end{align*}
	where the forth inequality is because of $p \le \theta$.
	By the definition of $\cC_{s,t}$, we have 
	\begin{align}
		\EE\left[ \norm{\tPi\tsb_t}^2 \right] \le \cC_{s,t}  
	\end{align}
	and
	\begin{align}
		\cC_{s,t+1} \le \left(1-\tilde{\theta}\right)\cC_{s,t} + 4\alpha p  \left(\norm{\nabla F(\xb_t) - \nabla F(\qb_t)}^2 + 2m\bsig^2\right).
	\end{align}
	
	Now, we are going to prove the consensus error related to $ \txb_t $.
	First, we have
	\begin{align*}
		&\EE\left[ \norm{\tPi\tW\tPi \left( \txb_t - \eta_t\Big( \tsb_t + \nabla F(\xb_t,\xi_t)_{\#} -\nabla F(\qb_t, \xi_\tau)_{\#}\Big)  \right)}^2 \right]\\
		\le& 
		\EE\left[\norm{ \tPi\tW\tPi (\txb_t - \eta_t \tsb_t) }^2\right] + \eta_t^2  \EE\left[\norm{ \tPi\tW\tPi \big( \nabla F(\xb_t, \xi_t)_{\#} -\nabla F(\qb_t, \xi_\tau)_{\#}\big)  }^2\right]\\
		& -2 \eta_t \EE\left[ \dotprod{\tPi\tW\tPi \txb_t,  \tPi\tW\tPi \big( \tsb_t +  \nabla F(\xb_t, \xi_t)_{\#} -\nabla F(\qb_t, \xi_\tau)_{\#}\big)} \right]
		\\
		=&
		\EE\left[\norm{ \tPi\tW\tPi (\txb_t - \eta_t \tsb_t) }^2\right] + \eta_t^2  \EE\left[\norm{ \tPi\tW\tPi \big(  \nabla F(\xb_t, \xi_t)_{\#} -\nabla F(\qb_t, \xi_\tau)_{\#}\big)  }^2\right]\\
		& -2  \eta_t\dotprod{\tPi\tW\tPi (\txb_t - \eta_t \tsb_t),  \tPi\tW\tPi \big(  \nabla F(\xb_t)_{\#} -\nabla F(\qb_t)_{\#}\big)} \\
		\le&
		\EE\left[\norm{ \tPi\tW\tPi (\txb_t - \eta_t \tsb_t) }^2\right] + \eta_t^2  \EE\left[\norm{ \tPi\tW\tPi \big(  \nabla F(\xb_t, \xi_t)_{\#} -\nabla F(\qb_t, \xi_\tau)_{\#}\big)  }^2\right]
		\\ &+ \frac{\tilde{\theta}}{2} \norm{ \tPi\tW\tPi (\txb_t - \eta_t \tsb_t) } + \frac{2\eta_t^2}{\tilde{\theta}} \norm{\tPi\tW\tPi \big(  \nabla F(\xb_t)_{\#} -\nabla F(\qb_t)_{\#}\big)}^2\\
		=&
		\left(1 + \frac{\tilde{\theta}}{2}\right) \norm{ \tPi\tW\tPi (\txb_t - \eta_t \tsb_t) }^2 +  \eta_t^2  \EE\left[\norm{ \tPi\tW\tPi \big(  \nabla F(\xb_t, \xi_t)_{\#} -\nabla F(\qb_t, \xi_\tau)_{\#}\big)  }^2\right]\\
		&+ \frac{2\eta_t^2}{\tilde{\theta}} \norm{\tPi\tW\tPi \big(  \nabla F(\xb_t)_{\#} -\nabla F(\qb_t)_{\#}\big)}^2\\
		\le&
		\left(1 + \frac{\tilde{\theta}}{2}\right) \left( \left(1 + \frac{\tilde{\theta}}{2}\right) \norm{\bPi \tW \bPi \txb_t}^2 + \left(1 + \frac{2}{\tilde{\theta}}\right) \eta_t^2 \norm{\bPi \tW \bPi \tsb_t}^2 \right)\\
		& +\eta_t^2  \EE\left[\norm{ \tPi\tW\tPi \big(  \nabla F(\xb_t, \xi_t)_{\#} -\nabla F(\qb_t, \xi_\tau)_{\#}\big)  }^2\right] + \frac{2\eta_t^2}{\tilde{\theta}} \norm{\tPi\tW\tPi \big(  \nabla F(\xb_t)_{\#} -\nabla F(\qb_t)_{\#}\big)}^2\\
		\le& 
		\left(1 + \frac{\tilde{\theta}}{2}\right)^2\norm{\bPi \tW \bPi \txb_t}^2 + \frac{6\eta_t^2}{\tilde{\theta}} \norm{\bPi \tW \bPi \tsb_t}^2\\
		& +\eta_t^2  \EE\left[\norm{ \tPi\tW\tPi \big(  \nabla F(\xb_t, \xi_t)_{\#} -\nabla F(\qb_t, \xi_\tau)_{\#}\big)  }^2\right] + \frac{2\eta_t^2}{\tilde{\theta}} \norm{\tPi\tW\tPi \big(  \nabla F(\xb_t)_{\#} -\nabla F(\qb_t)_{\#}\big)}^2\\
		=& 
		\left(1 + \frac{\tilde{\theta}}{2}\right)^2\EE\left[ \norm{\tPi\tW^2\tPi \left( \txb_{t-1} - \eta_{t-1}\Big( \tsb_{t-1} + \nabla F(\xb_{t-1},\xi_{t-1})_{\#} -\nabla F(\qb_{t-1}, \xi_\tau)_{\#}\Big)  \right)}^2 \right]\\
		&+ \frac{6\eta_t^2}{\tilde{\theta}} \norm{\bPi \tW \bPi \tsb_t}^2 +\eta_t^2  \EE\left[\norm{ \tPi\tW\tPi \big(  \nabla F(\xb_t, \xi_t)_{\#} -\nabla F(\qb_t, \xi_\tau)_{\#}\big)  }^2\right]\\
		&+\frac{2\eta_t^2}{\tilde{\theta}} \norm{\tPi\tW\tPi \big(  \nabla F(\xb_t)_{\#} -\nabla F(\qb_t)_{\#}\big)}^2,
	\end{align*}
	where the second and third inequality are because of Cauchy's inequality.
	
	Using above equation recursively, we can obtain that
	\begin{align*}
		&\EE\left[ \norm{\tPi\tW\tPi \left( \txb_t - \eta_t\Big( \tsb_t + \nabla F(\xb_t, \xi_t)_{\#} -\nabla F(\qb_t, \xi_\tau)_{\#}\Big)  \right)}^2 \right]\\
		\le& 
		\left(1 + \frac{\tilde{\theta}}{2}\right)^{2t} \norm{  \tPi\tW^{t+1}\tPi \txb_0}^2 + \frac{6}{\tilde{\theta}}\sum_{i=0}^{t} \eta_i^2 \left( 1 + \frac{\tilde{\theta}}{2}\right)^{2(t-i)} \norm{\tPi\tW^{t+1-i}\tPi \tsb_i}^2\\
		& +\sum_{i=0}^{t} \eta_i^2 \left( 1 + \frac{\tilde{\theta}}{2}\right)^{2(t-i)} \EE\left[\norm{\tPi\tW^{t+1-i}\tPi \Big( \nabla F(\xb_i, \xi_i)_{\#} -\nabla F(\qb_i,\xi_{\tau_i})_{\#} \Big)}^2\right]\\
		& + \frac{2}{\tilde{\theta}} \sum_{i=0}^t \eta_i^2 \left( 1 + \frac{\tilde{\theta}}{2}\right)^{2(t-i)} \norm{\tPi\tW^{t+1-i}\tPi \Big( \nabla F(\xb_i)_{\#} -\nabla F(\qb_i)_{\#} \Big)}^2\\
		\le&  
		\alpha\left(1 + \frac{\tilde{\theta}}{2}\right)^{2t} (1 - \tilde{\theta})^{2(t+1)} \norm{\tPi \txb_0}^2 + \frac{6\alpha}{\tilde{\theta}} \sum_{i=0}^{t} \eta_i^2 \left( 1 + \frac{\tilde{\theta}}{2}\right)^{2(t-i)} (1 - \tilde{\theta})^{2(t+1-i)} \norm{\tPi \tsb_i}^2\\
		&+\sum_{i=0}^{t} \eta_i^2 \left( 1 + \frac{\tilde{\theta}}{2}\right)^{2(t-i)} \left(1 -\tilde{\theta}\right)^{2(t+1-i)} \EE\left[\norm{\tPi \left( \nabla F(\xb_i, \xi_i)_{\#} -\nabla F(\qb_i, \xi_{\tau_i})_{\#} \right)}^2\right]\\
		&+\frac{2}{\tilde{\theta}} \sum_{i=0}^t \eta_i^2 \left( 1 + \frac{\tilde{\theta}}{2}\right)^{2(t-i)} \left(1 -\tilde{\theta}\right)^{2(t+1-i)}  \norm{\tPi \Big( \nabla F(\xb_i)_{\#} -\nabla F(\qb_i)_{\#} \Big)}^2\\
		\le& 
		\alpha\left(1-\frac{\tilde{\theta}}{2}\right)^{t+1}\norm{\tPi \txb_0}^2 + \frac{6\alpha}{\tilde{\theta}} \sum_{i=0}^{t} \eta_i^2  \left(1 - \frac{\tilde{\theta}}{2}\right)^{t+1-i} \norm{\tPi \tsb_i}^2\\
		&+\alpha\sum_{i=0}^{t} \eta_i^2 \left( 1 - \frac{\tilde{\theta}}{2}\right)^{t+1-i}  \EE\left[\norm{\tPi \left( \nabla F(\xb_i, \xi_i)_{\#} -\nabla F(\qb_i, \xi_{\tau_i})_{\#} \right)}^2\right]\\
		&+\frac{2\alpha}{\tilde{\theta}} \sum_{i=0}^t \eta_i^2 \left( 1 - \frac{\tilde{\theta}}{2}\right)^{t+1-i}  \norm{\tPi \Big( \nabla F(\xb_i)_{\#} -\nabla F(\qb_i)_{\#} \Big)}^2\\
		\stackrel{\eqref{eq:nn}}{\le}& 
		\alpha\left(1-\frac{\tilde{\theta}}{2}\right)^{t+1}\norm{\tPi \txb_0}^2 + \frac{6\alpha}{\tilde{\theta}} \sum_{i=0}^{t} \eta_i^2  \left(1 - \frac{\tilde{\theta}}{2}\right)^{t+1-i} \norm{\tPi \tsb_i}^2\\
		& +  \frac{3\alpha}{\tilde{\theta}} \sum_{i=0}^{t} \eta_i^2  \left(1 - \frac{\tilde{\theta}}{2}\right)^{t+1-i}  \norm{\tPi \Big( \nabla F(\xb_i)_{\#} -\nabla F(\qb_i)_{\#} \Big)}^2\\
		& + 2\alpha m\bsig^2 \sum_{i=0}^{t} \eta_i^2 \left( 1 - \frac{\tilde{\theta}}{2}\right)^{t+1-i} 
		\\
		=& \cC_{x,t+1}.
	\end{align*}
	By the definition of $\cC_{x,t+1}$, we can obtain the following inequality
	\begin{align*}
		&\cC_{x,t+1} \\
		\le& \left(1 - \frac{\tilde{\theta}}{2}\right) \cC_{x,t} +  \frac{3\alpha}{\tilde{\theta}} \left(1 - \frac{\tilde{\theta}}{2}\right) \eta_t^2 \left(  \norm{\tPi \tsb_t}^2 + 2 \norm{\nabla F(\xb_t) - \nabla F(\qb_t)}^2  \right) + 2\alpha m \bsig^2 \left(1 - \frac{\tilde{\theta}}{2}\right) \eta_t^2 \\
		\le& 
		\left(1 - \frac{\tilde{\theta}}{2}\right) \cC_{x,t} +  \frac{3\alpha}{\tilde{\theta}}  \eta_t^2 \left( \cC_{s,t} + 2 \norm{\nabla F(\xb_t) - \nabla F(\qb_t)}^2  \right) + 2\alpha m \bsig^2  \eta_t^2.
	\end{align*}
	By the update rule of $\txb_t$, it holds that
	\begin{align}
		\EE\left[\norm{\tPi\txb_{t+1}}^2\right] 
		=\EE\left[ \norm{\tPi\tW\tPi \left( \txb_t - \eta_t\big( \tsb_t + \zeta_t \left(\nabla F(\xb_t)_{\#} -\nabla F(\qb_t)_{\#}\right)\big)  \right)}^2 \right]
		\le \cC_{x,t+1}.
	\end{align}
	
\end{proof}

\subsection{Proof of Lemma~\ref{lem:psi_1}}
\begin{proof}[Proof of Lemma~\ref{lem:psi_1}]
	By the definition of $ \widetilde{\Psi}_{t+1} $, we have
\begin{align*}
	&\cC_{x,t+1} + C_{1, t+1} \cdot \cC_{s, t+1} + C_{2, t+1} \cdot \norm{\nabla F(\qb_{t+1}) - \nabla F(\mathbf{1}x^*)}^2\\
	\stackrel{\eqref{eq:ts}\eqref{eq:tx}\eqref{eq:qx}}{\le}&
	\left(1 -\frac{\tilde{\theta}}{2} \right) \cC_{x,t} + \left(1 - \frac{\tilde{\theta}}{2} + \frac{3\alpha\eta_t^2}{\tilde{\theta} C_{1, t}}\right) C_{1,t} \cdot \cC_{s, t}
	+ (1 -p) C_{2,t} \norm{\nabla F(\qb_t) - \nabla F(\mathbf{1}x^*)}^2\\
	& + \left(\frac{6\alpha\eta_t^2}{\tilde{\theta}} + 4\alpha p C_{1,t}\right) \norm{\nabla F(\xb_t) - \nabla F(\qb_t)}^2 + \left(2 \alpha m \eta_t^2 + 8\alpha m p C_{1,t} \right)\bsig^2\\
	& + 4mLpC_{2,t} \Big(f(\bbx_t) - f(x^*)\Big) + 2pL^2 C_{2,t} \norm{\bPi \xb_t }^2\\
	\stackrel{\eqref{eq:xq}\eqref{eq:tsx}}{\le}&
	\left(1 - \frac{\tilde{\theta}}{2} + \frac{24\alpha L^2 \eta_t^2}{\tilde{\theta}} + 16\alpha p L^2 C_{1,t} + 2pL^2 C_{2,t}\right) \cdot \cC_{x,t}\\
	& +\left(1 - \frac{\tilde{\theta}}{2} + \frac{3\alpha\eta_t^2}{\tilde{\theta} C_{1, t}}\right) C_{1,t} \cdot \cC_{s, t}\\
	& + \left(1 - p + \frac{12\alpha\eta_t^2}{\tilde{\theta} C_{2,t}} + 8\alpha p \cdot \frac{C_{1,t}}{C_{2,t}} \right) C_{2,t} \norm{\nabla F(\qb_t) - \nabla F(\mathbf{1}x^*)}^2\\
	& + \left(2 \alpha m \eta_t^2 + 8\alpha m p C_{1,t} \right)\bsig^2 + \left(4mLpC_{2,t} + \frac{48\alpha m L\eta_t^2 }{\tilde{\theta}} + 32\alpha p mL C_{1,t} \right)\Big(f(\bbx_t) - f(x^*)\Big)\\
	\le&
	\left(1 - \frac{\tilde{\theta}}{4}\right) \left( \cC_{x,t} +  C_{1,t} \cdot \cC_{s, t} + C_{2,t} \norm{\nabla F(\qb_t) - \nabla F(\mathbf{1}x^*)}^2 \right)\\
	& \frac{2^{12}\cdot 3^2\cdot m \eta_t^2}{\tilde{\theta}} \cdot \bsig^2 + \frac{2^{12}\cdot 3^2\cdot m L \eta_t^2}{\tilde{\theta}} \cdot \Big(f(\bbx_t) - f(x^*)\Big),
\end{align*}
where the last inequality is due to the setting of parameters.
\end{proof}

\subsection{Proof of Lemma~\ref{lem:main_tld}}
\begin{proof}[Proof of Lemma~\ref{lem:main_tld}]
By Lemma~\ref{lem:equal}, we can conclude that Lemma~\ref{lem:main_2} still holds. 
\begin{align*}
&\EE\left[ \norm{\bbx_{t+1} - x^*}^2 + \frac{48L\eta_{t+1}}{m\tilde{\theta}} \widetilde{\Psi}_{t+1} \right]\\
\stackrel{\eqref{eq:xx_a}\eqref{eq:psi_1}}{\le}&
\left(1 - \frac{\mu\eta_t}{2}\right) \norm{\bbx_t - x^*}^2 - 2\eta_t \left(1 - 2\eta_t L \right)\big(f(\bbx_t) - f(x^*)\big) \\
&+ \eta_t^2\cdot \frac{\bsig^2}{m}
+\frac{2L\eta_t \left( 1 + 2\eta_t L \right) }{m} \norm{\bPi \xb_t}^2 
+ \left(1- \frac{\tilde{\theta}}{4}\right) \frac{48L\eta_t}{m\tilde{\theta}} \widetilde{\Psi}_t \\
& + \frac{2^{16}\cdot 3^2 \cdot  L\eta_t^3}{\tilde{\theta}^2}\bsig^2 + \frac{2^{16}\cdot 3^2 \cdot  L^2\eta_t^3}{\tilde{\theta}^2} \cdot \big(f(\bbx_t) - f(x^*)\big)\\
\le&
\left(1 - \frac{\mu\eta_t}{2}\right) \norm{\bbx_t - x^*}^2 - 2\eta_t \left(1 - 2\eta_t L - \frac{2^{15}\cdot 3^2 \cdot L^2\eta_t^2}{\tilde{\theta}^2}\right)\big(f(\bbx_t) - f(x^*)\big)\\
&+\left(1- \frac{\tilde{\theta}}{4}\right) \frac{48L\eta_t}{m\tilde{\theta}} \widetilde{\Psi}_t + \frac{6L\eta_t \cdot \cC_{x,t}}{m} + \eta_t^2\cdot \frac{\bsig^2}{m} + \frac{2^{16}\cdot 3^2 \cdot  L\eta_t^3}{\tilde{\theta}^2} \cdot \bsig^2\\
\stackrel{\eta_t \le \frac{\tilde{\theta}}{2^8\cdot 3\cdot L}}{\le}& 
\left(1 - \frac{\mu\eta_t}{2}\right) \left(\norm{\bbx_t - x^*}^2 + \frac{48L\eta_t}{m\tilde{\theta}} \cdot \widetilde{\Psi}_t\right) - \frac{7\eta_t}{8} \Big(f(\bbx_t) - f(x^*)\Big) + \eta_t^2\cdot \frac{\bsig^2}{m} + \frac{2^{16}\cdot 3^2 \cdot  L\eta_t^3}{\tilde{\theta}^2} \cdot \bsig^2\\
\le&
\exp\left(-\frac{\mu\eta_t}{2}\right)\left(\norm{\bbx_t - x^*}^2 + \frac{48L\eta_t}{m\tilde{\theta}} \cdot \widetilde{\Psi}_t\right) - \frac{7\eta_t}{8} \Big(f(\bbx_t) - f(x^*)\Big) + \eta_t^2\cdot \frac{\bsig^2}{m} + \frac{2^{16}\cdot 3^2 \cdot  L\eta_t^3}{\tilde{\theta}^2} \cdot \bsig^2,
\end{align*}
where  the last inequality is because of $ 1 - x \le \exp(-x) $ for $0\le x <1$.
\end{proof}

\section{Decentralized Stochastic Gradient Tracking}

\begin{algorithm}[tb]
	\caption{Decentralized Stochastic Gradient Tracking}
	\label{alg:dsgt}
	\begin{algorithmic}
		\STATE {\bfseries{Input}}: $x_0$, mixing matrix $W$, initial step size $\eta$.
		\STATE {\bfseries{Initialization}:} Set $\xb_0 = \mathbf{1}x_0$, $\qb_0 = \mathbf{1}x_0$, $\bs_0^{(i)} = \nabla f_i(\xb_0^{(i)},\xi_0)$, in parallel for $i \in [m]$, $\tau = 0$.
		\FOR {$t = 1,\dots, T$}
		\STATE Update
		\begin{align}
			\xb_{t+1} =& W \left( \xb_t - \eta_t  \bs_t \right) \label{eq:dsgt_up}, \\
			\bs_{t+1} =& W\bs_t +  \nabla F(\xb_{t+1},\xi_{t+1}) - \nabla F(\xb_t,\xi_t). \label{eq:ss_GT}
		\end{align}
		\ENDFOR
	\end{algorithmic}
\end{algorithm}

\end{appendix}

\end{document}